\input ./style/arxiv-general.cfg
\documentclass[MSNbibl,citesort,number,seceqn,dvips]{arxbj}
\makeatletter
   \@ifpackageloaded{graphicx}{}{\usepackage{graphicx}}
\makeatother
\usepackage{upgreek}

\aid{0}
\volume{21}
\issue{3}
\pubyear{2015}
\firstpage{1435}
\lastpage{1466}
\doi{10.3150/14-BEJ610} % kopijuoti is 'New paper accepted'
\makeatletter

\newcommand{\rright}{\right}
\newcommand{\lleft}{\left}

\newtheorem{Lem}{Lemma}[section]
\newtheorem{Theor}{Theorem}[section]
\newcommand{\I}{I}
\newcommand{\E}{\operatorname{E}}
\newcommand{\R}{\mathbb{R}}

\newcommand{\s}{\mathcal{S}}
\newcommand{\ny}{n\rightarrow\infty}
\newcommand{\thetab}{\boldsymbol{\theta}}
\newcommand{\thetabv}{\boldsymbol{\vartheta}}
\newcommand{\taub}{\boldsymbol{\tau}}

\newcommand{\argmin}{\operatorname{\arg\min}}
\newcommand{\tr}{\operatorname{tr}}
\newcommand{\diag}{\operatorname{diag}}

\newcommand{\XX}{{\boldsymbol{\mathcal{X}}}}
\newcommand{\FZ}{\mathrm{FZ}}

\newcommand{\Var}{\operatorname{Var}}
\renewcommand{\vec}{\operatorname{vec}}

\newproclaim{Assumption}{Assumption}
\newremark{Remark}{Remark}[section]
\makeatother

\begin{document}
\begin{frontmatter}

\title{Local bilinear multiple-output quantile/depth regression}
\runtitle{Local multiple-output quantile regression}

\begin{aug}
%%%% inicialai - be tarpu
\author[1,2]{\inits{M.}\fnms{Marc}~\snm{Hallin}\corref{}\thanksref{1,2,e1}\ead[label=e1,mark]{mhallin@ulb.ac.be}},
\author[3]{\inits{Z.}\fnms{Zudi} \snm{Lu}\thanksref{3}\ead[label=e2]{zudilu@gmail.com}\ead[label=u2,url]{https://sites.google.com/site/zudiluwebsite/}},
\author[1]{\inits{D.}\fnms{Davy} \snm{Paindaveine}\thanksref{1,e3}\ead[label=e3,mark]{dpaindav@ulb.ac.be}\ead
[label=u3,url]{http://homepages.ulb.ac.be/\textasciitilde dpaindav}} 
\and
\author[4]{\inits{M.}\fnms{Miroslav} \snm{\v{S}iman}\thanksref{4}\ead[label=e4]{siman@utia.cas.cz}}
%%\runauthor{} %% auto
%\dedicated{}
\address[1]{E.C.A.R.E.S., Universit\'{e} libre de Bruxelles
CP114/04,
50 Avenue F.D. Roosevelt,
B-1050 Brussels, Belgium.\\
\printead{e1,e3,u3}}

\address[2]{ORFE, Princeton University,
Sherrerd Hall,
Princeton, NJ 08544,
USA}

\address[3]{Mathematical Sciences \& Statistical Sciences Research Institute,
University of Southampton,
Building 54, Highfield,
Southampton SO17 1BJ, UK.
\printead{e2};\\\printead{u2}}

\address[4]{Department of Stochastic Informatics, Institute of Information Theory and Automation of the ASCR,
Pod Vod\'{a}renskou v\v{e}\v{z}\'{i} 4,
CZ-182 08 Prague 8, Czech Republic.
\printead{e4}}
\end{aug}

% HISTORY:
\received{\smonth{8} \syear{2012}}
\revised{\smonth{5} \syear{2013}}

% ABSTRACT
%
\begin{abstract}
A new quantile regression concept, based on a directional version of
Koenker and Bassett's traditional single-output one, has been
introduced in [\textit{Ann. Statist.} (2010) \textbf{38} 635--669] for multiple-output location/linear
regression problems. The polyhedral contours provided by the empirical
counterpart of that concept, however, cannot adapt to unknown nonlinear
and/or heteroskedastic dependencies. This paper therefore introduces
local constant and local linear (actually, bilinear) versions of those
contours, which both allow to asymptotically recover the conditional
halfspace depth contours that completely characterize the response's
conditional distributions. Bahadur representation and asymptotic
normality results are established. Illustrations are provided both on
simulated and real data.
\end{abstract}

% KEYWORDS
% visi is mazosios raides ir pagal abecele
%

\end{frontmatter}

%s1 #&#
\section{Introduction}

%s1.1 #&#
\subsection{Quantile/depth contours: From multivariate location to
multiple-output regression}
\label{sec1.1}

A multiple-output extension of Koenker and Bassett's celebrated concept
of regression quantiles was recently proposed in Hallin, Paindaveine,
and \v{S}iman \cite{Hetal10} (hereafter HP\v{S}). That extension
provides regions that are enjoying, at population level, a double
interpretation in terms of quantile and halfspace depth regions. In the
empirical case, those regions are limited by polyhedral contours which
can be computed via parametric linear programming techniques.

Those results establish a strong and quite fruitful link between two
seemingly unrelated statistical worlds -- on one hand the typically
one-dimensional concept of quantiles, deeply rooted into the strong
ordering features of the real line and L$_1$ optimality, with linear
programming algorithms, and traditional central-limit asymptotics; the
intrinsically multivariate
concept of depth on the other hand, with geometric characterizations,
computationally intensive combinatorial algorithms, and nonstandard asymptotics.
From their relation to depth, quantile hyperplanes and regions inherit
a variety of geometric properties -- connectedness, nestedness,
convexity, affine-equivariance \ldots  while, via its relation to
quantiles, depth accedes to L$_1$ optimality, feasible linear
programming algorithms, and tractable asymptotics.

The HP\v{S} approach, however, is focused on the case of i.i.d.
$m$-variate observations $\mathbf{Y}_1,\ldots, \mathbf{Y}_n$, and
the quantile/depth contours they propose provide a consistent
reconstruction of the corresponding population contours in $\mathbb
{R}^m$ -- call them \textit{unconditional} or \textit{location}
contours. In the presence of covariates $\mathbf{X}_1,\ldots, \mathbf
{X}_n$, with $\mathbf{X}_i=(1,\mathbf{W}_i^{\prime})^{\prime}$, the
objective
of the statistical analysis is a study of the influence of the
covariate(s) $\mathbf{W}$ on the response $\mathbf{Y}$, that is, a
study of the distribution of $\mathbf{Y}$ conditional on $\mathbf
{W}$. The contours of interest, thus, are the collection of the
population \textit{conditional quantile/depth contours} of $\mathbf
{Y}$, indexed by the values $\mathbf{w}\in\R^{p-1}$ of $\mathbf
{W}$ -- that is, the collection of location ($p=1$) quantile/depth
contours associated with the conditional
(on $\mathbf{W}= \mathbf{w}$) distributions of $\mathbf{Y}$.

An apparently simple solution would consist in introducing the
covariate values $\mathbf{w}$ into the linear equations that
characterize (via the minimization of an L$_1$ criterion) the HP\v{S}
contours. The resulting regions and contours, unfortunately, in general
carry little information about conditional distributions, and rather
produce some averaged (over the covariate space) quantile/depth
contours -- the only exception being the overly restrictive case of a
linear regression relation between the response and the covariates,
under which, for some ${\mathbf{\frak{b}}}\in\mathbb{R}^p$, the
distribution of $\mathbf{Y} -(\mathbf{1},\mathbf{w}^{\prime}){\mathbf{\frak
{b}}} $
conditional on $\mathbf{W} = \mathbf{w}$ does not depend on $\mathbf
{w}\in\R^{p-1}$.

This problem is not specific to the multiple-output context and, in the
traditional single-output setting, it has motivated weighted, local
polynomial and nearest-neighbor versions of quantile regression, among
others. We refer to \cite{YJ97,YJ98,YL04} for conceptual insight and
practical information, to \cite{CX09,EV09,HS09,Hetal09,Ketal09,ZW09}
for some recent asymptotic results,
% in relevant contexts,
and to \cite{BG90,C91,Getal03,GS05,H00,H04,I04,WS94} for some less
recent ones.

Our objective in this paper is to extend those local estimation ideas
to the HP\v{S} concept of multiple-output regression quantiles. Since
local constant and local linear methods have been shown to perform
extremely well in the single-output single-regressor case (Yu and Jones
\cite{YJ98}), we will concentrate on local constant and local \textit
{bilinear} approaches -- in the multiple-output context, indeed, it
turns out that the adequate extensions of locally linear procedures are
of a \textit{bilinear} nature. Just as in the single-output case, the
local methods we propose in this paper do not require any a priori
knowledge of any trend and -- see \cite{KZ10} for details -- asymptotically characterize the conditional distributions
of $\mathbf{Y}$ given $\mathbf{W}= \mathbf{w}$ for any $\mathbf
{w}\in\mathbb{R}^{p-1}$. The final result is thus much more
informative on the dependence of $\mathbf{Y}$ on the covariates than
any standard linear or local polynomial mean regression. %This is
%illustrated in Figures \ref{figcone}(c)-(d), which shows a very good
%agreement with the population quantities provided in Figure
%\ref{figcone}(e); in particular, both trend and heteroskedascticity
%components are now appropriately recovered.

It should be clear, however, that our methods, as well as other local
nonparametric methods, do not escape the curse of dimensionality, and
will run into problems in the presence of high-dimensional regressors.
It follows indeed from the asymptotic results of Section~\ref{asympsec} and, more particularly, from the rates in Theorem~\ref
{a.normality}, that consistency rates are affected by $p$ but not by $m$.

Growth chart applications (with $(p-1)=1$) do not suffer this drawback,
as only univariate kernels %in the local constant approach,
%$m+1)$-dimensional ones in the local bilinear approach,
are involved. Growth charts (reference curves, percentile curves) have
been used for a long time by practitioners in order to assess the
impact of regressors on the quantiles of some given univariate variable
of interest, and several methods have been developed (see, e.g., \cite{BR96,C88,WH06,WR97}, and the references therein),
including single-response quantile regression (see \cite
{Getal02,Wetal06}). Much less results are available in the
multiple-output case, with a recent proposal by Wei \cite{W08}, who
defines a new concept of dynamic multiple-output regression contours
generalizing single-output proposals by \cite{CX09}, \cite{K07} and
\cite{WH06}. These contours, however, do not have the nature and
interpretation of (conditional) depth contours. They enjoy interesting
conditional coverage probability properties (without any ``minimal
volume'' or ``maximal accuracy'' features, though) but rely on a
sequential conditioning of response components, and crucially depend on
the order adopted for that conditioning. Their empirical versions are
equivariant under marginal location-scale transformations of the
response, but they are neither affine- nor rotation-equivariant. Our
methodology, which is based on entirely different principles, appears
as a natural alternative (see \cite{Metal11} for a real-data example
of bivariate growth charts based on the methods we are describing
here), yielding affine-equivariant regression contours with
well-accepted conditional depth interpretation; moreover, % ontrary to
%\cite{W08},
we provide consistency and asymptotic distributional results.

%s1.2 #&#
\subsection{Motivating examples}
% \subsection{Simulated data}
\label{simul1}

\begin{figure}

\includegraphics{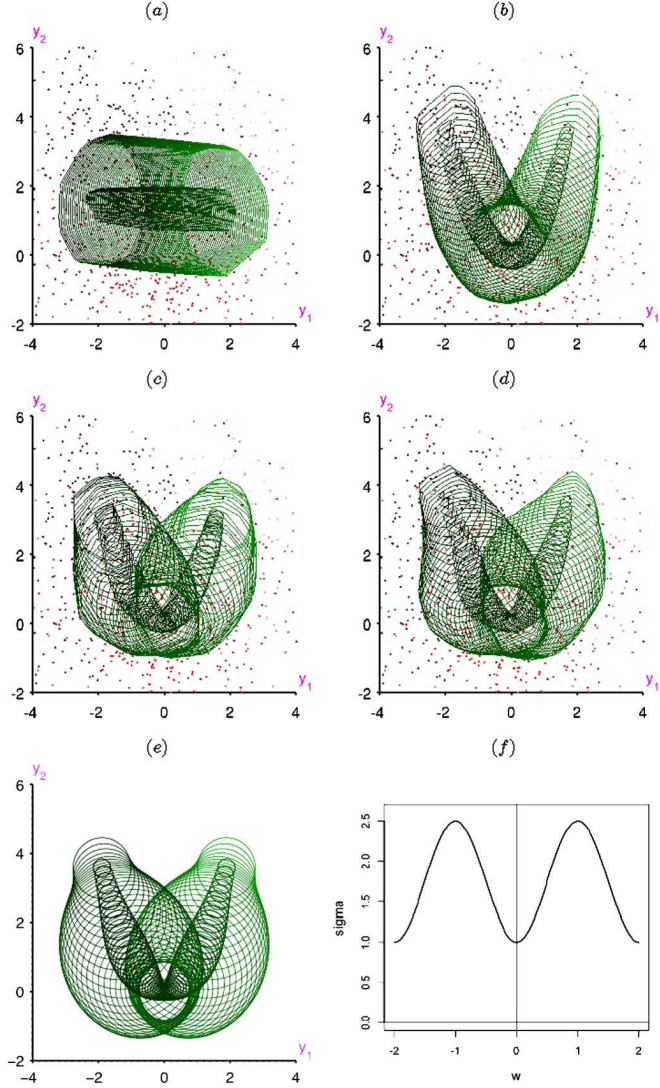}

\caption{For $n=999$ points following the model $(Y_1,Y_2) = (W,W^2) +
(1+\frac
{3}{2}(\sin(\frac{\uppi}{2}W))^2){\boldsymbol{\varepsilon}}$, where
$W\sim
U([-2,2])$ and ${\boldsymbol{\varepsilon}}\sim\mathcal{N}(0,1)^2$ are
independent, the plots above show the intersections, with hyperplanes
orthogonal to the $w$-axis at fixed $w$-values $-1.89, -1.83,-1.77,
\dots, 1.89$, of
(a)
the HP\v{S} regression quantile regions with the single random regressor
$W$,
(b)
the HP\v{S} regression quantile regions
with random regressors
$W$ and $W^2$,
and
(c)--(d)
the proposed local constant and local bilinear regression quantile
regions (in each case, $\tau=0.2$ and $\tau=0.4$ are considered). For
the sake of comparison, the corresponding population (conditional)
halfspace depth regions are provided in (e). The conditional scale
function $w\mapsto1+\frac{3}{2}(\sin(\frac{\uppi}{2}w))^2$ is
plotted in (f). Local methods use a Gaussian kernel and bandwidth value
$H=0.37$, and 360 equispaced directions $\mathbf{u}\in\mathcal{S}^1$
were used to obtain results in (d).}
\label{figcone}
\end{figure}

As a motivating example, we generated $n=999$ points from the model
\[
(Y_1,Y_2) = \bigl(W,W^2 \bigr) + \biggl(1+
\frac{3}{2} \biggl(\sin \biggl(\frac{\uppi}{2}W \biggr)
\biggr)^2 \biggr){\boldsymbol {\varepsilon}},
\]
with $W\sim U([-2,2])$ independent of the bivariate standard normal
vector ${\boldsymbol{\varepsilon}}$. In Figure~\ref{figcone}, we are plotting
the $\tau=0.2$ and $\tau=0.4$ HP\v{S} regression quantile/depth
contours obtained by using the covariate vector $\mathbf{X}=(1,W)'$
(Figure~\ref{figcone}(a)) and the covariate vector $\mathbf
{X}=(1,W,W^2)'$ (Figure~\ref{figcone}(b)) in the equations of the
quantile/depth hyperplanes of the (global) HP\v{S} mehod. More
precisely, these figures provide the intersections of the HP\v{S}
contours with hyperplanes orthogonal to the $w$-axis at fixed
$w$-values $-1.89, -1.83,-1.77, \dots, 1.89$.
%
%f1 #&#
%
%f2 #&#
\begin{figure}

\includegraphics{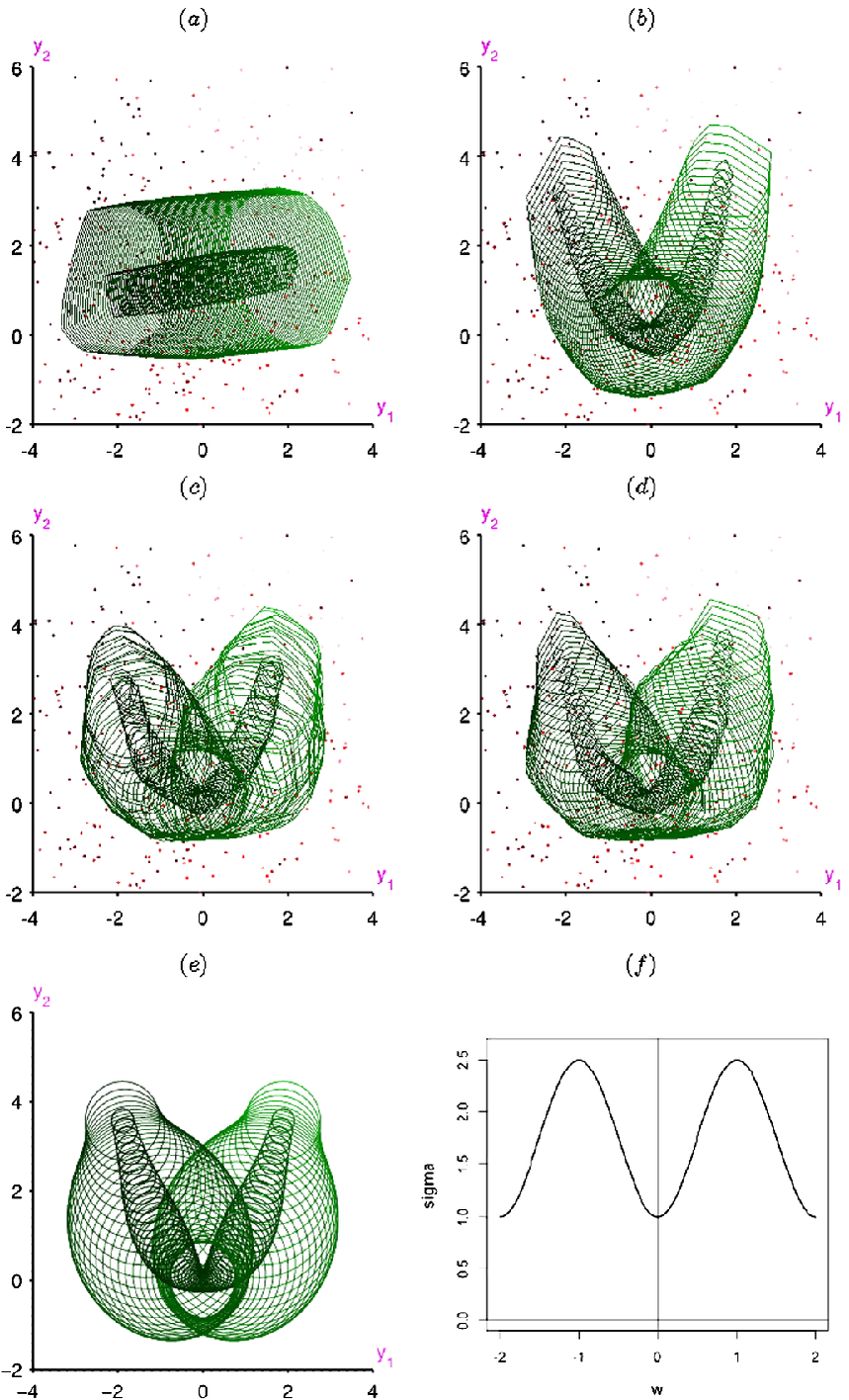}

\caption{For $n=499$ points following the model $(Y_1,Y_2) = (W,W^2) +
(1+\frac
{3}{2}(\sin(\frac{\uppi}{2}W))^2){\boldsymbol{\varepsilon}}$, where
$W\sim
U([-2,2])$ and ${\boldsymbol{\varepsilon}}\sim\mathcal{N}(0,1)^2$ are
independent, the plots above show the intersections, with hyperplanes
orthogonal to the $w$-axis at fixed $w$-values $-1.89, -1.83,-1.77,
\dots, 1.89$, of
(a)
the HP\v{S} regression quantile regions with the single random regressor
$W$,
(b)
the HP\v{S} regression quantile regions
with random regressors
$W$ and $W^2$,
and
(c)--(d)
the proposed local constant and local bilinear regression quantile
regions (in each case, $\tau=0.2$ and $\tau=0.4$ are considered). For
the sake of comparison, the corresponding population (conditional)
halfspace depth regions are provided in (e). The conditional scale
function $w\mapsto1+\frac{3}{2}(\sin(\frac{\uppi}{2}w))^2$ is
plotted in (f). Local methods use a Gaussian kernel and bandwidth value
$H=0.37$, and 360 equispaced directions $\mathbf{u}\in\mathcal{S}^1$
were used to obtain results in (d).}
\label{figcone2}
\end{figure}

Clearly, the results are very poor: Figure~\ref{figcone}(a) neither
reveals the parabolic trend, nor the periodic heteroskedasticity
pattern in the data. Although it is obtained by fitting the ``true''
regression function, Figure~\ref{figcone}(b), while doing much better
with the trend, still fails to catch heteroskedasticity correctly.
Instead of providing genuine conditional quantile/depth contours, the
``global'' HP\v{S} methodology produces some averaged (over the $w$
values) contours.

In contrast, the contours obtained from the local constant and local
bilinear methods proposed in this paper -- without exploiting any a
priori knowledge of the actual regression function -- exhibit a very
good agreement with the population contours (see Figure~\ref{figcone}(c)--(e) to which we refer for details); both the parabolic
trend and the periodic heteroskedascticity features now are picked up
quite satisfactorily. Note that, compared to the local constant
approach, the local bilinear one does better, as expected, close to the
boundary of the regressor space (in particular, the local constant
approach is missing the decay of the conditional scale when $w$
converges to $-2$).

\begin{figure}

\includegraphics{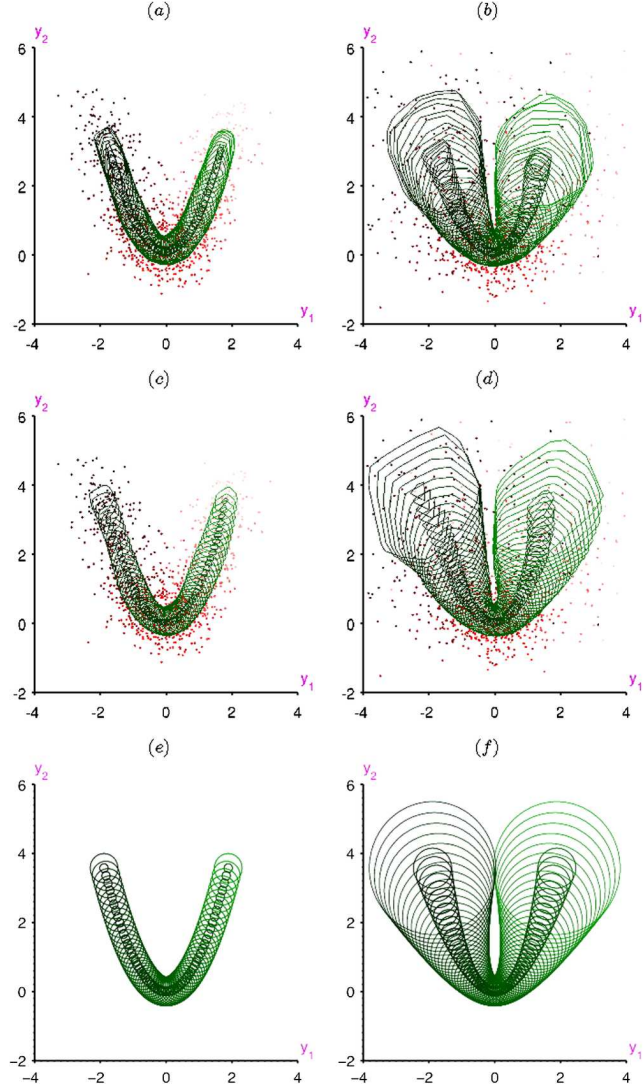}

\caption{Local multiple-output quantile regression with Gaussian
kernel and ad-hoc bandwidth $H = 0.37$: cuts through $w \in\{
-1.89, -1.83, -1.77, \dots, 1.89\}$ for $\tau= 0.2$ and $\tau= 0.4$
corresponding to $n = 999$ random points drawn from a homoskedastic
model $(Y_1,Y_2) = (W,W^2) + {\boldsymbol{\varepsilon}}$ ((a), (c))
or a
heteroskedastic model $(Y_1,Y_2) = (W,W^2) + (1+W^2){\boldsymbol
{\varepsilon}
}$ ((b), (d)), where $W\sim U([-2,2])$ and $ {\boldsymbol{\varepsilon
}} \sim
\mathcal{N}(0,1/4)^2$ are independent. The plots are showing the
intersections, with hyperplanes orthogonal to the $w$-axis at fixed
$w$-values, of the contours obtained either from the local constant
method ((a), (b)) or the local bilinear one ((c), (d)). Color scaling
of the points (resp., the intersections) mimics their regressor values,
whose higher values are indicated by lighter red (resp., lighter
green). For the sake of comparison, the population (conditional)
halfspace depth regions are provided in (e) and (f).
A color version of this figure is more readable, and can be found in the on-line edition of the paper.}

\label{DemoFig}
\end{figure}

Similar comments remain valid for smaller sample sizes; see Figure~\ref{DemoFig}, based on a sample of $n=499$ data points.

A second example is contrasting a homoskedastic setup and a
heteroskedastic one. More specifically, we generated $n = 999$ points
from the homoskedastic model
$(Y_1,Y_2) = (W,W^2) + {\boldsymbol{\varepsilon}}$ and from the
heteroskedastic one $(Y_1,Y_2) = (W,W^2) + (1+W^2){\boldsymbol
{\varepsilon}}$,
where $W\sim U([-2,2])$ and $ {\boldsymbol{\varepsilon}} \sim
\mathcal
{N}(0,1/4)^2$ are mutually independent. As above, the intersections of
the resulting contours with hyperplanes orthogonal to the $w$-axis at
fixed $w$-values are provided. Figure~\ref{DemoFig} shows those
intersections for the local constant and local bilinear quantile
contours associated with $w \in\{-1.89, -1.83, -1.77, \dots, 1.89\}$,
for $\tau= 0.2$ and $\tau= 0.4$. As in the previous example, those
sample contours approximate their population counterparts (shown in Figure~\ref{DemoFig}(e) and (f)) remarkably well. In particular, the inner
regions mimic the trend faithfully even for quite extreme regressor
values. Again, the local bilinear method seems to provide a much better
boundary behavior than its local constant counterpart; in the heteroskedastic case, the latter
indeed severely underestimates the conditional scale for extreme values
of $W$.
%
%%%%%%%%%%%%%%%%%%%%%%%%%%%%%%%%%%%%%%%%%%
%
%f3 #&#

%s1.3 #&#
\subsection{Relation to the depth and multivariate quantile literature
%%, and application to multiple-output growthcharts
}
\label{13}

As already explained, this work is lying at the intersection of two
distinct, if not unrelated, strands of the statistical
literature -- namely (i) statistical depth and (ii) multivariate
quantiles. Under both strands, definitions have been proposed for
\textit{unconditional} concepts, that is, for statistical models that
do not involve covariates. When covariates are present, the focus is
shifted from unconditional features to \textit{conditional} ones. The
main objective, indeed, now is the analysis of the dependence of a
response $\mathbf{Y}$ on a set of covariates $\mathbf{X}$, that is, a
study of the distributions of $\mathbf{Y}$ conditional on the values
$\mathbf{x}$ of $\mathbf{X}$ --  in its broadest sense, the \textit
{regression problem} -- and various attempts have been made %, in both
%strands of literature,
to propose \textit{regression} versions of (unconditional) depth or
quantile concepts, respectively.

Now, if a study of the dependence on $\mathbf{x}$ of the distributions
of $\mathbf{Y}$ conditional on $\mathbf{X}=\mathbf{x}$ is the main
objective, \textit{conditional} depth and \textit{conditional}
(multivariate) quantiles, associated with the distributions of $\mathbf
{Y}$ conditional on $\mathbf{X}=\mathbf{x}$, are or should be the
concepts of interest. Not all definitions of \textit{regression depth}
or {(multiple-output) regression quantiles} are meeting that
requirement, though. Nor do they all preserve, conditionally on
$\mathbf{X}=\mathbf{x}$, the distinctive properties of a
depth/quantile concept. In contrast with this, the concept we are
proposing in this paper, being the conditional version of the
unconditional HP\v{S} concept, enjoys all the properties that are
expected from a conditional depth/quantile concept, while fully
characterizing the conditional distributions of~$\mathbf{Y}$.
%Its empirical counterpart moreover is compuconsistently estimates its
%conditional
%
% Our main contribution is the definition, in the multiple-output
%\textit{regression} context, of a depth/quantile concept that fully
%characterizes the dependence (the conditional distributions) of the
%multivariate response on covariates, while preserving the main
%properties that are expected from a depth/quantile concept.

%s1.3.1 #&#
\subsubsection{Regression depth}\label{subsecdep} %Much attention has
%been given, in the depth-related literature, to the relation between
%depth and quantile ideas (the converse is not true).
An excellent summary of depth-related problems is provided in Serfling
\cite{S10}, which further clarifies the nature of depth by placing it in
the broader perspective of the so-called \textit{DOQR paradigm},
relating Depth to the companion concepts of Outlyingness, Quantiles,
and Ranks. To the best of our knowledge, this paradigm never has been
considered in a conditional (regression) context, but it seems quite
desirable that any \textit{regression depth} concept should similarly
be placed, conditionally, in the same DOQR perspective.

The celebrated \emph{regression depth} concept by Rousseeuw and Hubert
\cite{RH99}, %(see \citet{RH99} and \citet{RR99}),
for instance, does not bear any direct relation to conditional depth
and the DOQR paradigm. Rather than the depth of a point in the
observation space, that concept aims at defining, via \textit
{non-fits} and \textit{breakdown values}, the depth of a
(single-output) regression hyperplane. A multiple-output version is
considered in Bern and Eppstein \cite{BD02}. Similarly, an elegant
general theory has
been developed by Mizera \cite{Miz02} who, in the context of a general
parametric model, defines the depth of a parameter value. Again, the
approach and, despite the terminology, the concept, is of a different
nature, unrelated to any conditional depth. Extensions to a
nonparametric regression setting, moreover, seem problematic.

Kong and Mizera \cite{KM10} propose an approach to unconditional depth,
based on projection quantiles, which provides an approximation to the
unconditional halfspace depth contours -- see \cite{Hetal10} and \cite{KM10}. Although an
application to bivariate growth charts is briefly
described, in which a local smoothing, based on regression spline
techniques, of their unconditional concept is performed (little details
are provided), the regression setting is only briefly touched there. In
particular, no asymptotic analysis of the type we are providing in
Section~\ref{asympsec} is made available.

%s1.3.2 #&#
\subsubsection{Multivariate regression quantiles}\label{subsecquant}

Turning to conditional multivariate or multiple-output regression
quantile issues, much work has been devoted to the notion of \emph
{spatial} regression quantiles; see, essentially, Chakraborty
\cite{Ch03} for
linear and Cheng and De~Gooijer \cite{CG07} for nonparametric
regression. %Those
%multiple-output spatial regression quantiles nevertheless
Despite a strong depth flavor, those spatial quantiles and spatial
regression quantiles, however, intrinsically fail to be
affine-equivariant; Chakraborty \cite{Ch03} defines
affine-equivariant spatial
quantiles for linear regression via a transformation--retransformation
device, but, to the best of our knowledge, there exists no
affine-equivariant version of spatial quantiles for general
nonparametric regression.

%It is then appealing to try and extend the HP\v S location
%multivariate quantiles to the regression context. The linear
%regression extension that was already defined in HP\v S, however, is
%not satisfactory, as it was shown through the motivating examples
%above. The goal of this paper is therefore to improve on this by
%proposing a local extension of the HP\v S concept, that addresses
%successfully the general nonparametric regression problem.

For the sake of completeness, one also should mention here the closely
related literature on growth charts described at the end of Section~\ref{sec1.1}, which, besides a lack of affine-invariance, essentially
fails, in the multiple-output case, to address the conditional nature
of the regression quantile concept it is dealing with.
\subsection{Outline of the paper}The rest of this paper is organized
as follows.
Section~\ref{definsec} defines the (population) conditional regression
quantile/depth regions and contours we would like to estimate in the sequel.
This estimation will make use of (empirical) \textit{weighted}
multiple-output regression quantiles, which we introduce in Section~\ref{weightsec}.
Section~\ref{locsec} explains how these weighted quantiles lead to
local constant (Section~\ref{locconstsec}) and local bilinear (Section~\ref{loclinsec}) depth contours.
Section~\ref{asympsec} provides asymptotic results (Bahadur
representation and asymptotic normality) both for the local constant
and local bilinear cases. Section~\ref{Bandsec} deals with the
practical problem of bandwidth selection.
In Section~\ref{pracsec}, the usefulness and applicability of the
proposed methods are illustrated on real data.
%Section~\ref{finsec}
%
Finally, the \hyperref[app]{Appendix} collects proofs of asymptotic results.

%%%%%%%%%%%%%%%%%%%%%%%%%%%%%%%%%%%%%%%%%%

%s2 #&#
\section{Conditional multiple-output quantile/depth contours}
\label{definsec}
%%%%%%%%%%%%%%%%%%%%%%%%%%%%%%%%%%%%%%%%%%

Denote by $(\mathbf{X}_i^{\prime}, \mathbf{Y}_i^{\prime})^{\prime
}= (X_{i1} ,
\ldots, X_{ip}$, $ Y_{i1},\ldots, Y_{im})^{\prime}$, $i=1,\ldots,
n$, an
observed $n$-tuple of independent copies of $(\mathbf{X}^{\prime},
\mathbf
{Y}^{\prime})^{\prime}$, where $\mathbf{Y}:=(Y_1,\dots,Y_m)^{\prime
}$ is an
$m$-dimensional response and $\mathbf{X}:=(1,\mathbf{W}^{\prime
})^{\prime}$ a
$p$-dimensional random vector of covariates. For any $\tau\in(0,1)$
and any direction $\mathbf{u}$ in the unit sphere $\mathcal{S}^{m-1}$
of the $m $-dimensional space of the response $\mathbf{Y}$, the HP\v
{S} concept produces a hyperplane $\pi_{\tau\mathbf{u}}$ ($\pi^{(n)}
_{\tau\mathbf{u}}$ in the empirical case) which is defined as the
classical Koenker and Bassett regression quantile hyperplane of order
$\tau$ once $(\mathbf{0}_{p-1}^{\prime},\mathbf{u}^{\prime
})^{\prime}$ has been
chosen as the ``vertical direction'' in the computation of the relevant
L$_1$ deviations.

More specifically, decompose $\mathbf{y}\in\R^m$ into $(\mathbf
{u}^{\prime}\mathbf{y})\mathbf{u} + {\boldsymbol{\Gamma}}_\mathbf
{u}({\boldsymbol{\Gamma}}_\mathbf{u}^{\prime}\mathbf{y})$, where
${\boldsymbol{\Gamma}
}_\mathbf{u}$ is such that $ (\mathbf{u},{\boldsymbol{\Gamma
}}_\mathbf
{u} )$ is an $m\times m$ orthogonal matrix; then the \textit
{directional quantile hyperplanes} $\pi_{\tau\mathbf{u}}$ and $\pi
^{(n)}_{\tau\mathbf{u}}$ are the hyperplanes with equations
%
%e2.1 #&#
\begin{equation}
\label{ddprems} \mathbf{u}^{\prime}\mathbf{y} - \mathbf{c}_{\boldsymbol{\tau}}
^{\prime}{\boldsymbol{\Gamma} }_\mathbf{u}^{\prime}\mathbf{y}-
\mathbf{a}^{\prime}_{\boldsymbol
{\tau}} \bigl(1,\mathbf {w}^{\prime}
\bigr)^{\prime}=0 %$$
%and
%$$
\quad \mbox{and}\quad
\mathbf{u}^{\prime}\mathbf{y} - \mathbf{c}^{(n)\prime
}_{\boldsymbol{\tau}} {
\boldsymbol{\Gamma}}_\mathbf{u}^{\prime}\mathbf{y}-\mathbf
{a}^{(n)\prime}_{\boldsymbol{\tau}} \bigl(1,\mathbf{w}^{\prime}
\bigr)^{\prime}=0
\end{equation}
($\mathbf{w}\in\R^{p-1}$) minimizing, with respect to $\mathbf
{c}\in\R^{m-1}$ and $\mathbf{a}\in\R^{p}$,
%
%e2.2 #&#
\begin{equation}
\label{unweightedfirst} \mathrm{E} \bigl[\rho_\tau \bigl(\mathbf{u}^{\prime}
\mathbf{Y} - \mathbf {c}^{\prime} {\boldsymbol{\Gamma}}_\mathbf{u}
^{\prime}\mathbf{Y}-\mathbf {a}^{\prime}\mathbf{X} \bigr) \bigr]
%$$
\quad \mbox{and}\quad %$$
\sum_{i=1}^n
\rho_\tau \bigl(\mathbf{u}^{\prime}\mathbf{Y}_i -
\mathbf {c}^{\prime} {\boldsymbol{\Gamma}}_\mathbf{u}^{\prime}
\mathbf{Y}_i-\mathbf {a}^{\prime}\mathbf{X}_i
\bigr),
\end{equation}
respectively,
where $\zeta\mapsto\rho_\tau(\zeta) $, with
%
%e2.3 #&#
\begin{equation}
\label{check} \rho_\tau(\zeta) := \zeta \bigl(\tau- \I[\zeta<0]
\bigr) = \max \bigl\{( \tau -1)\zeta,\tau\zeta \bigr\} = \bigl(|\zeta|+(2\tau-1)
\zeta \bigr)/2 ,\qquad \zeta\in\R
\end{equation}
as usual denotes the
well-known $\tau$-quantile \textit{check function}. HP\v{S} moreover
show that $\pi_{\tau\mathbf{u}}$ and $\pi^{(n)}_{\tau\mathbf{u}}$
can equivalently be defined, in a more symmetric way, as the
hyperplanes with equations
%
%e2.4 #&#
\begin{equation}
\label{dd} \mathbf{b}^{\prime}_{\boldsymbol{\tau}} \mathbf{y}-\mathbf
{a}^{\prime}_{\boldsymbol{\tau}} \bigl(1,\mathbf{w}^{\prime}
\bigr)^{\prime}=0 \quad \mbox{and}\quad \mathbf {b}^{(n)\prime}_{\boldsymbol{\tau}}
\mathbf{y}-\mathbf{a}^{(n)\prime}_{\boldsymbol{\tau}} \bigl(1,\mathbf
{w}^{\prime} \bigr)^{\prime}=0, % \mathbf{y}\in\R^m, \mathbf{w}\in\R^{p-1}
\end{equation}
minimizing, with respect to $\mathbf{b}\in\R^{m}$ satisfying
$\mathbf{b}^{\prime}\mathbf{u}=1$ and $\mathbf{a}\in\R^{p}$, the
L$_1$ criteria
%
%e2.5 #&#
\begin{equation}
\label{unweighted} \mathrm{E} \bigl[\rho_\tau \bigl(\mathbf{b}^{\prime}
\mathbf{Y}-\mathbf {a}^{\prime} \mathbf{X} \bigr) \bigr] %$$
\quad
\mbox{and}\quad %$$
\sum_{i=1}^n
\rho_\tau \bigl(\mathbf{b}^{\prime}\mathbf{Y}_i-
\mathbf {a}^{\prime} \mathbf{X}_i \bigr),
\end{equation}
respectively.

For $p=1$, the multiple-output regression model reduces to a
multivariate location one: $\mathbf{a} _{\boldsymbol{\tau}}$ and
$\mathbf
{a}^{(n)}_{\boldsymbol{\tau}}$ reduce to scalars, ${ a} _{\boldsymbol
{\tau}}$ and ${
a}^{(n)}_{\boldsymbol{\tau}}$, while the equations describing
$\pi_{\tau\mathbf{u}}$ and $\pi^{(n)}_{\tau\mathbf{u}}$ take the
simpler forms
%
%e2.6 #&#
\begin{equation}
\label{locpi} \mathbf{u}^{\prime}\mathbf{y} - \mathbf{c}_{\boldsymbol{\tau}}
^{\prime}{\boldsymbol{\Gamma} }_\mathbf{u}^{\prime}
\mathbf{y}-{ a}_{\boldsymbol{\tau}} =0 \quad \mbox{and}\quad \mathbf{u}^{\prime}
\mathbf{y} - \mathbf{c}_{\boldsymbol{\tau
}}^{(n)\prime}{\boldsymbol{
\Gamma}}_\mathbf{u}^{\prime}\mathbf {y}-{ a}^{(n)}_{\boldsymbol{\tau}}
=0,
\end{equation}
respectively. Those \textit{location} quantile hyperplanes $\pi
_{\tau\mathbf{u}}$ and $\pi^{(n)}_{\tau\mathbf{u}}$ are studied in
detail in HP\v{S}, where it is shown that their fixed-$\tau$
collections characterize regions and contours that actually coincide
with the Tukey halfspace depth ones. Consistency, asymptotic normality
and Bahadur-type representation results for the $\pi^{(n)}_{\tau
\mathbf
{u}}$'s are also provided there, together with a linear programming
method for their computation.

The objective here is an analysis of the distribution of $\mathbf{Y}$
conditional on $\mathbf{W}$, that is, of the dependence of $\mathbf
{Y}$ on $\mathbf{W}$ -- in strong contrast with traditional regression,
where investigation is limited to the mean of $\mathbf{Y}$ conditional
on $\mathbf{W}$. The relevant quantile hyperplanes, depth regions and
contours of interest are the location quantile/depth
hyperplanes/regions/contours associated (in the sense of HP\v{S}) with
the $m$-dimensional distributions of $\mathbf{Y}$ conditional on
$\mathbf{W}$ -- more precisely, with the distributions $\mathrm
{P}^{\mathbf{Y}\vert\mathbf{W}=\mathbf{w}_0}$ of $\mathbf{Y}$
conditional on $\mathbf{W}=\mathbf{w}_0$ ($\mathbf{w}_0\in\R
^{p-1}$). We now carefully define these objects -- call them % %that we
%will call
\emph{$\mathbf{w}_0$-conditional $\boldsymbol{\tau}$-quantile or depth
hyperplanes, regions and contours}.
%
%The problem (when $p\geq2$, and under absolutely continuous $\bf W$)
%is that the empirical counterparts of those $\mathbf{w}_0$-conditional
%concepts are degenerate (since the hyperplane with equation $
%\mathbf{w}=\mathbf{w}_0$ is of measure zero in $\R^{p-1}$), and this
%is why a \textit{local} approach is required.

%\subsection{Conditional multiple-output quantiles and depth}
%\label{conditsec}

Let $\tau\in(0,1)$ and $\mathbf{u}\in\mathcal{S}^{m-1} := \{
\mathbf{u} \in\R^m\dvtx  \|\mathbf{u}\| = 1\}$ (the unit sphere in $\R
^m$), and write ${\boldsymbol{\tau}} := \tau\mathbf{u}$. Denoting by
$\mathbf{w}_0 $ some fixed point of $\R^{p-1}$ at which the marginal
density $f^\mathbf{W}$ of $\mathbf{W}$ does not vanish (in order for
the distribution of $\mathbf{Y}$ conditional on $\mathbf{W}=\mathbf
{w}_0$ to make sense), define the \textit{extended} and \textit
{restricted $\mathbf{w}_0$-conditional $\boldsymbol{\tau}$-quantile
hyperplanes of $\mathbf{Y}$} as the $(m+p-2)$-dimensional and
$(m-1)$-dimensional hyperplanes
%
%e2.7 #&#
\begin{equation}
\label{condquantext} {\boldsymbol{\pi}}_{{\boldsymbol{\tau}};\mathbf{w}_0} := \bigl\{ \bigl(
\mathbf{w}^{\prime}, \mathbf {y}^{\prime} \bigr)^{\prime}\in
\R^{p-1}\times\R^m \mid \mathbf{b}_{{\boldsymbol{\tau}};\mathbf{w}_0}^{\prime}
\mathbf{y} - a_{{\boldsymbol{\tau}
};\mathbf{w}_0} = 0 \bigr\}
\end{equation}
and
%
%e2.8 #&#
\begin{equation}
\label{condquantrestr} \pi_{{\boldsymbol{\tau}};\mathbf{w}_0} := \bigl\{ \bigl(\mathbf {w}_0^{\prime},
\mathbf{y}^{\prime} \bigr)^{\prime}\in\R^{p-1}\times\R
^m \mid \mathbf{b}_{{\boldsymbol{\tau}};\mathbf{w}_0}^{\prime}\mathbf{y} -
a_{{\boldsymbol{\tau}};\mathbf{w}_0} = 0 \bigr\},
\end{equation}
respectively, where $a_{{\boldsymbol{\tau}};\mathbf{w}_0}$ and
$\mathbf{b}_{{\boldsymbol{\tau}};\mathbf{w}_0}$
minimize
%
%e2.9 #&#
\begin{equation}
\label{CondOpt} \Psi_{\tau;\mathbf{w}_0}(a,\mathbf{b}) : = \E \bigl[
\rho_\tau \bigl(\mathbf {b}^{\prime}\mathbf{Y}- a \bigr) \mid
\mathbf{W}=\mathbf{w}_0 \bigr] \qquad \mbox{subject to }
\mathbf{b}^{\prime}\mathbf{u}=1,
\end{equation}
with the check function $\rho_\tau$ defined in (\ref{check}).
Comparing (\ref{CondOpt}) with (\ref{unweighted}) immediately shows that
$\pi_{{\boldsymbol{\tau}};\mathbf{w}_0}$ is the $(m-1)$-dimensional
(location) $\boldsymbol{\tau}$-quantile hyperplane of $\mathbf{Y}$ associated
with the distribution of $\mathbf{Y}$ conditional on $\mathbf
{W}=\mathbf{w}_0$. Of course, $\pi_{{\boldsymbol{\tau}};\mathbf
{w}_0}$ is
also the intersection of ${\boldsymbol{\pi}}_{{\boldsymbol{\tau
}};\mathbf{w}_0}$
with the $m$-dimensional hyperplane $C_{\mathbf{w}_0}:=\{(\mathbf
{w}_0^{\prime}, \mathbf{y}^{\prime})^{\prime}\mid\mathbf{y}\in
\R^m \}$. This,
and the fact that ${\boldsymbol{\pi}}_{{\boldsymbol{\tau}};\mathbf
{w}_0}$ is
``parallel to the space of
covariates'' \label{horizpage} (in the sense that if $(\mathbf
{w}_0^{\prime}
,\mathbf{y}_0^{\prime})^{\prime}\in{\boldsymbol{\pi
}}_{{\boldsymbol{\tau}};\mathbf{w}_0}$,
then $(\mathbf{w}^{\prime},\mathbf{y}_0^{\prime})^{\prime}\in
{\boldsymbol{\pi}}_{{\boldsymbol{\tau}};\mathbf{w}_0}$ for all~$\mathbf{w}$), fully characterizes
${\boldsymbol{\pi}}_{{\boldsymbol{\tau}};\mathbf{w}_0}$.

Associated with ${\boldsymbol{\pi}}_{{\boldsymbol{\tau}};\mathbf
{w}_0}$ are the
extended \textit{upper and lower $\mathbf{w}_0$-conditional
${\boldsymbol{\tau}}$-quantile halfspaces}
\[
\mathbf{ H}^{+}_{{\boldsymbol{\tau}};\mathbf{w}_0} := \bigl\{ \bigl(\mathbf
{w}^{\prime},\mathbf {y}^{\prime} \bigr)^{\prime}\in
\R^{p-1}\times\R^m \mid\mathbf {b}_{{\boldsymbol{\tau}
};\mathbf{w}_0}^{\prime}
\mathbf{y} - a_{{\boldsymbol{\tau
}};\mathbf{w}_0} \geq 0 \bigr\}
\]
and
\[
\mathbf{ H}^{-}_{{\boldsymbol{\tau}};\mathbf{w}_0} := \bigl\{ \bigl(\mathbf
{w}^{\prime},\mathbf {y}^{\prime} \bigr)^{\prime}\in
\R^{p-1}\times\R^m \mid\mathbf {b}_{{\boldsymbol{\tau}
};\mathbf{w}_0}^{\prime}
\mathbf{y} - a_{{\boldsymbol{\tau
}};\mathbf{w}_0} < 0 \bigr\},
\]
with the extended (cylindrical) \textit{$\mathbf{w}_0$-conditional
quantile/depth regions}
%
%e2.10 #&#
\begin{equation}
\label{2.2} \mathbf{ R}_{\mathbf{w}_0}(\tau) := \bigcap
_{\mathbf{u}\in\s^{m-1}} \bigl\{ \mathbf{ H}^{+}_{\tau\mathbf{u};\mathbf{w}_0}
\bigr\}
\end{equation}
and their boundaries $ \partial\mathbf{ R}_{\mathbf{w}_0}(\tau)$, the
extended \textit{$\mathbf{w}_0$-conditional quantile/depth contours}.
The intersections of those extended
regions $\mathbf{ R}_{\mathbf{w}_0}(\tau)$ (resp., contours $
\partial
\mathbf{ R}_{\mathbf{w}_0}(\tau)$) with $C_{\mathbf{w}_0}$ are the
restricted $\mathbf{w}_0$-conditional quantile/depth regions
${R}_{\mathbf{w}_0}(\tau)$ (resp., contours $\partial{R}_{\mathbf
{w}_0}(\tau)$), that is, the location HP\v{S} regions (resp.,
contours) for $\mathbf{Y}$, conditional on $\mathbf{W}=\mathbf
{w}_0$. It follows from HP\v{S} that those regions are compact,
convex, and nested. As a consequence, the regions $\mathbf{
R}_{\mathbf
{w}_0}(\tau)$ also are closed, convex, and nested.

Finally, define the \textit{nonparametric $\tau$-quantile/depth
regions} as
\[
\mathbf{ R}( \tau) := \bigcup _{\mathbf{w}_0\in\R^{p-1}} R_{\mathbf
{w}_0}(
\tau) = \bigcup _{\mathbf{w}_0\in\R^{p-1}} \bigl(\mathbf{ R}_{\mathbf{w}_0}(
\tau)\cap C_{\mathbf{w}_0} \bigr)
\]
and write $\partial\mathbf{ R}(\tau)$ for their boundaries. The regions
$\mathbf{ R}( \tau)$ are still closed and nested but they adapt to the
general dependence of $\mathbf{Y}$ on $\mathbf{W}$: in particular,
$\partial\mathbf{ R}(\tau)$, for any $\tau$, goes through \textit{all}
corresponding $\partial{R}_{\mathbf{w}_0}(\tau)$'s, $\mathbf
{w}_0\in\R^{p-1}$. Consequently, the regions $\mathbf{ R}( \tau)$ in
general are no longer
convex.

The fixed-$\mathbf{w}_0$ collection (over $\tau\in(0, 1/2 ))$ of
all $\mathbf{w}_0$-conditional location quantile/depth contours
$\partial{R}_{\mathbf{w}_0}(\tau)$ (which, by construction, are the
intersections of $\partial\mathbf{ R}( \tau)$ with the ``vertical
hyperplanes'' $C_{\mathbf{w}_0}$) will be called a \textit{$\mathbf
{w}_0$-quantile/depth cut} or \textit{$\mathbf{w}_0$-cut}.
Such cuts are of crucial interest, since they characterize the
distribution of $\mathbf{Y}$ conditional on $\mathbf{W}=\mathbf
{w}_0$, hence provide a full description of the dependence of the
response $\mathbf{Y}$ on the regressors $\mathbf{W}$. Note that the
nonparametric contours $ \partial\mathbf{ R}(\tau)$, via the location
depth interpretation, for fixed $\mathbf{w}_0$, of the $\partial
{R}_{\mathbf{w}_0}(\tau)$'s, inherit a most interesting
interpretation as ``regression depth contours''. Clearly, this concept
of regression depth, that defines regression depth of any point
$(\mathbf{w}',\mathbf{y}')'\in\mathbb{R}^{m+p-1}$, is not of the
same nature as the regression depth concept proposed in \cite{RH99},
that defines the depth of any regression ``fit'' (i.e., of any
regression hyperplane).

%%%%%%%%%%%%%%%%%%%%%%%%%%%%%%%%%%%%%%%%%%

%s3 #&#
\section{Weighted multiple-output empirical quantile regression}
\label{weightsec}

Under the assumption of absolute continuity,
the number of observations, in a sample of size $n$, belonging to
$C_{\mathbf{w}_0}$ clearly is (a.s.) zero, which implies that no
empirical version of the conditional regression hyperplanes (\ref
{condquantext}) or (\ref{condquantrestr}) can be constructed. If
nonparametric $\tau$-quantile/depth regions or contours, or simply
some selected cuts, are to be estimated, local smoothing techniques
have to be considered. Those techniques typically involve weighted
versions, with sequences ${\boldsymbol{\omega}}^{(n)}_{\mathbf
{w}_0}=(\omega^{(n)}
_{\mathbf{w}_0,i},i=1,\ldots,n)$ of weights, of the empirical
quantile regression hyperplanes developed in HP\v{S}. In this section,
we provide general definitions and basic results for such weighted
concepts, under fixed sample size $n$ and weights $\omega_i$; see
\cite{Hluetal09} for another approach combining weights with
halfspace depth. In Section~\ref{locsec}, we will consider the
data-driven weights to be used in the local approach.

Consider a sample of size $n$, with observations $(\mathbf
{X}_i^{\prime},
\mathbf{Y}_i^{\prime})^{\prime}= ((1,\mathbf{W}_i^{\prime}) ,
\mathbf{Y}_i^{\prime})^{\prime}
$, $i=1,\ldots,n$, along with $n$ nonnegative weights ${\omega}_i$
satisfying (without any loss of generality) $\sum_{i=1}^n {\omega}_i
= n$ ($\omega_i\equiv1$ then yields the unweighted case). The
definitions of HP\v{S}
extend, \textit{mutatis mutandis}, quite straightforwardly, into the
following weighted versions.
The coefficients ${\mathbf{a}_{{\boldsymbol{\tau}};{\omega
}}^{(n)}}\in\R
^{p}$ and $ {\mathbf{b}_{{\boldsymbol{\tau}};\omega}^{(n)} } \in\R
^{m}$ of
the \textit{weighted empirical $\boldsymbol{\tau}$-quantile hyperplane}
%
%e3.1 #&#
\begin{equation}
\label{vectora} {\boldsymbol{\pi}}^{(n)}_{{\boldsymbol{\tau}};\omega} := \bigl\{
\bigl( \mathbf {w}^{\prime}, \mathbf {y}^{\prime}
\bigr)^{\prime}\in \R^{p-1}\times\R^m \mid{\mathbf
{b}_{{\boldsymbol{\tau}
};{\omega}}^{(n)\prime} } \mathbf{y}-{\mathbf{a}_{{\boldsymbol
{\tau}
};{\omega}}^{(n)\prime}}
\bigl(1, \mathbf{w}^{\prime} \bigr)^{\prime}= 0 \bigr\}
\end{equation}
(an $(m+p-2)$-dimensional hyperplane) are defined as
the minimizers of
%
%e3.2 #&#
\begin{equation}
\label{weightobj} \Psi_{\tau;{\omega}}^{(n)}(\mathbf{a},\mathbf{b}) :=
\frac{1}{n}\sum_{i=1}^n {
\omega}_i\rho_\tau \bigl(\mathbf{b}^{\prime
}\mathbf
{Y}_i-\mathbf{a}^{\prime}\mathbf{X}_i \bigr)
\qquad \mbox{subject to } \mathbf{b}^{\prime}\mathbf{u}=1.
\end{equation}
As usual in the empirical case, the solution may not be unique, but the
minimizers always form a convex set.
When substituted for the ${\boldsymbol{\pi}}_{{\boldsymbol{\tau
}};\mathbf{w}_0}$'s
in the definitions of upper and lower conditional ${\boldsymbol{\tau}
}$-quantile halfspaces, those ${\boldsymbol{\pi}}^{(n)}_{{\boldsymbol
{\tau}};\omega}$'s
also characterize upper and lower weighted ${\boldsymbol{\tau
}}$-quantile halfspaces
$
\mathbf{ H}^{(n) +}_{{\boldsymbol{\tau}};\omega}$ and $
\mathbf{ H}^{(n) -}_{{\boldsymbol{\tau}};\omega}$, with \textit
{weighted $\tau
$-quantile/depth regions and contours}
\[
\mathbf{ R}^{(n)}_{\omega}(\tau) := \bigcap
_{\mathbf{u}\in\s
^{m-1}} \bigl\{ \mathbf{ H}^{(n) +}_{\tau\mathbf{u};\omega}
\bigr\} \quad \mbox{and}\quad \partial\mathbf{ R}^{(n)}_{\omega}(
\tau),
\]
respectively.
Note that the objective function in (\ref{weightobj}) rewrites as
\[
\Psi_{\tau;{\omega}}^{(n)}(\mathbf{a},\mathbf{b}) = \frac{1}{n}
\sum_{i=1}^n \rho_\tau \bigl(
\mathbf{b}^{\prime}\mathbf {Y}_{i;{\omega}}-\mathbf{a}^{\prime}
\mathbf{X}_{i;{\omega}} \bigr) ,
\]
with $\mathbf{X}_{i;{\omega}}:={\omega}_i\mathbf{X}_{i}$ and
$\mathbf{Y}_{i;{\omega}}:={\omega}_i\mathbf{Y}_{i}$. As an
important consequence, the weighted quantile/depth hyperplanes,
contours and regions can be computed in the same way as their
non-weighted counterparts because the corresponding algorithm in \cite
{PS11b} allows to have $(\mathbf{X}_i)_1\neq1$. Due to %the
quantile crossing, %phenomenon,
however, and contrary to the population regions and contours defined in
the previous section, the $\mathbf{ R}^{(n)}_{\omega}(\tau)$'s need
not be
nested for $p\geq2$; if nestedness is required, one may rather
consider the regions $\mathbf{ R}^{(n)}_{{\omega}\cap}(\tau) :=
\bigcap_{0<t\leq\tau} \{ \mathbf{ R}^{(n)}_{\omega}(t)\}$.

The necessary sample subgradient conditions for $({\mathbf
{a}_{{\boldsymbol{\tau}};\omega}^{(n)\prime}} ,{\mathbf
{b}_{{\boldsymbol{\tau}};{\omega
}}^{(n)\prime}} )^{\prime}$ can be derived as in the unweighted case.
They state in particular that
\begin{eqnarray*}
%\label{SampleProb}
\frac{1}{n}\sum_{i=1}^n
{\omega}_i \I \bigl[ {\mathbf{b}_{{\boldsymbol{\tau}};\omega}^{(n)\prime}
\mathbf{Y}}_i- \mathbf{a}_{{\boldsymbol{\tau}};\omega}^{(n)\prime}
\mathbf{X}_i < 0 \bigr] \leq\tau\leq\frac{1}{n}\sum
_{i=1}^n {\omega}_i \I \bigl[ {
\mathbf{b}_{{\boldsymbol{\tau}};\omega}^{(n)\prime} \mathbf{Y}}_i-
\mathbf{a}_{{\boldsymbol{\tau}};\omega}^{(n)\prime} \mathbf{X}_i \leq0 \bigr],
\end{eqnarray*}
which controls the probability contents of $\mathbf{
H}^{(n)-}_{{\boldsymbol{\tau}
};{\omega}}$ with respect to the distribution putting probability mass
$\omega_i/n$ on $(\mathbf{W}_i^{\prime}, \mathbf{Y}_i^{\prime
})^{\prime}$,
$i=1,\ldots,n$.
The width of this interval %in (\ref{D ^{(T)}_{t;m})
depends only on the weights ${\omega}_i$ associated with those data
points $(\mathbf{W}_i ^{\prime}, \mathbf{Y}_i^{\prime})^{\prime}$
that belong to
${\boldsymbol{\pi}}^{(n)}_{\taub;{\omega}}$.
Another consequence worth mentioning is that there always exists a
${\boldsymbol{\pi}}^{(n)}_{\tau\mathbf{u};{\omega}}$ hyperplane
containing at
least $(m+p-1)$ data points of the form $(\mathbf{W}_{i},\mathbf
{Y}_{i})$. With probability one, thus, the intersection defining
the regions $\mathbf{ R}^{(n)}_{\omega}(\tau)$ is finite.

Note that, unlike the extended conditional quantile hyperplanes (\ref
{condquantext}), the weighted empirical quantile hyperplanes (\ref
{vectora}) involve an unrestricted coefficient $\mathbf{a}\in\R
^{p}$. As a consequence, ${\boldsymbol{\pi}}^{(n)}_{\taub;{\omega}}
$ is not
necessarily parallel to the space of
covariates (as defined in page \pageref{horizpage}). That degree of
freedom will be exploited in the {local linear} approach described in
Section~\ref{loclinsec} (in an augmented regressor space, though,
which makes it bilinear rather than linear). If we impose the
additional constraint $\mathbf{a}=(a_1,0,\ldots, 0)^{\prime}$ in
(\ref
{vectora}) and (\ref{weightobj}), we obtain hyperplanes of the form
%
%e3.3 #&#
\begin{equation}
\label{scalara} {\boldsymbol{\pi}}^{(n)}_{{\boldsymbol{\tau}};\omega} := \bigl\{
\bigl( \mathbf {w}^{\prime}, \mathbf {y}^{\prime}
\bigr)^{\prime}\in \R^{p-1}\times\R^m \mid {
\mathbf{b}_{{\boldsymbol{\tau}};{\omega}}^{(n)\prime} } \mathbf {y}-{{a}_{1;{\boldsymbol{\tau}};\omega}^{(n)}}
= 0 \bigr\}.
\end{equation}
The corresponding minimization problem yields hyperplanes that are
parallel to the space of
covariates, hence ``horizontal'' cylindrical regions and contours, %.
%For the sake of simplicity, we avoid introduce any specific notation
%for them; such cylindrical contours
to be considered in the {local constant} approach of Section~\ref{locconstsec}.

Finally, it should be pointed out that ($\mathbf{y}$ and/or $\mathbf
{w}$)-affine-invariant weights $\omega_i:=\omega(\mathbf
{w}_i,\mathbf{y}_i)$ yield weighted quantile/depth hyperplanes,
regions, and contours with good ($\mathbf{y}$ and/or $\mathbf
{w}$)-affine-equivariance properties.

%%%%%%%%%%%%%%%%%%%%%%%%%%%%%%%%%%%%%%%%%%
%s4 #&#
\section{Local quantile/depth regression}
\label{locsec}
%%%%%%%%%%%%%%%%%%%%%%%%%%%%%%%%%%%%%%%%%%

%s4.1 #&#
\subsection{From weighted to local quantile/depth regression}
\label{weightdefsec}

The weighted quantiles of Section~\ref{weightsec} have an interest on
their own. They can be used for handling multiple identical
observations (allowing, for instance, for bootstrap procedures), or for
downweighting observations that are suspected to be outliers or
leverage points. Above all, weighted regression quantiles
allow for a nonparametric approach to regression quantiles that will
take care of the drawbacks of the unweighted approach of HP\v{S} (see
the example considered in Section~\ref{simul1}). In particular,
adequate sequences of weights will allow to estimate the conditional
contours described in Section~\ref{definsec}, thus extending to the
multiple-output case the \textit{local constant} and \textit{local
linear} approaches to quantile regression proposed, for example, by
\cite{YJ97,YJ98} in the single-output context.

The basic idea is very standard: in order to estimate $\mathbf
{w}_0$-conditional quantile/depth hyperplanes, regions or contours, we
will consider weighted quantile/depth hyperplanes, regions or contours,
with sequences of weights $\omega^{(n)}_i:=\omega^{(n)}_{\mathbf
{w}_0}(\mathbf{W}_i)$ based on weight functions of the form
%
%e4.1 #&#
\begin{equation}
\label{wdef} \mathbf{w}\mapsto{\omega}^{(n)}_{\mathbf{w}_0}(
\mathbf{w}) := %
%\det(
%\mathbf{H}\n_0)^{-1}
h_n^{-p+1}K
\bigl(h_n^{-1}(\mathbf{w}-\mathbf{w}_0)
\bigr),
\end{equation}
where $h_n$ %$\mathbf{H}\n_0$
is a sequence of %symmetric
positive %definite $(p-1)\times(p-1)$
bandwidths % matrices
and $K$ %is a
a nonnegative kernel function over $\R^{p-1} $.
%, that is, satisfies $\int_{\R^{p-1}} K(\mathbf{w})\mathrm{d}
%\mathbf{w}=1$.
The literature proposes a variety of possible kernels, and there is no
compelling
reason for not considering the most usual, such as the rectangular
(uniform), Epanechnikov or (spherical) Gaussian ones.
%\begin{enumerate}%{align*}
%\item[(i)] the rectangular (uniform) kernel $\DS{
%K_1(\mathbf{w}) = 2^{-(p-1)}\I[\mathbf{w} \in[-1,1]^{p-1}] ,}$
%\item[(ii)] the Epanechnikov % (parabolic)
% kernel $\DS{
%K_2(\mathbf{w}) = \frac{(p^2-1)\Gamma(\frac{p-1}{2})}{4\pi^{(p-1)/2}}
%(1-\mathbf{w}'\mathbf{w}) \I[\mathbf{w}'\mathbf{w} \leq1] ,}$ or
%\item[(iii)] the (spherical) Gaussian kernel $\DS{
%K_3(\mathbf{w}) = (2\pi)^{-(p-1)/2} \exp({-\mathbf{w}'\mathbf{w}}/2)
%.}$
%\end{enumerate}%{align*}
%As for the bandwidth matrices, we will restrict to the simple scalar
%case, that is, to $\mathbf{H}_0\n=h_n\mathbf{I}_{p-1}$ and write $K_h(
%\mathbf{w} - \mathbf{w}_0)$ for the weight ${\omega}\n_{\mathbf{w}_0}(
%\mathbf{w})$.

Since we typically intend, for any fixed $\tau\in(0,1)$, to compute
by means of parametric programming the directional quantile hyperplanes
for all $\mathbf{u} \in\mathcal{S}^{m-1}$, we should use the same
weights for all of them. This is why we only consider $\mathbf
{u}$-independent %(actually, even ${\boldsymbol{\tau}}$-independent)
bandwidths. However, exact computation of all quantiles (for each fixed
$\tau$) is possible in the local constant case, but not in the local
bilinear one. In the latter case, depth contours will be approximated
by sampling the unit sphere (in Figures \ref{figcone} and \ref{figcone2}, for instance, 360 directions were sampled uniformly over the unit circle), which of course would allow $\mathbf
{u}$-dependent bandwidths if desired.

%The weights considered above cover both kernel and nearest-neighbor
%quantile regression but exclude more sophisticated techniques such as
%double-kernel-, supersmoother- or LOWESS-based modifications. On the
%other hand, the choice of weights has no impact on computational
%issues, and special kernels (and bandwidths) can be selected for
%extreme $\mathbf{w}_0$'s to take care of boundary effects, for
%instance.

%s4.2 #&#
\subsection{Local constant quantile/depth contours}
\label{locconstsec}

\begin{figure}

\includegraphics{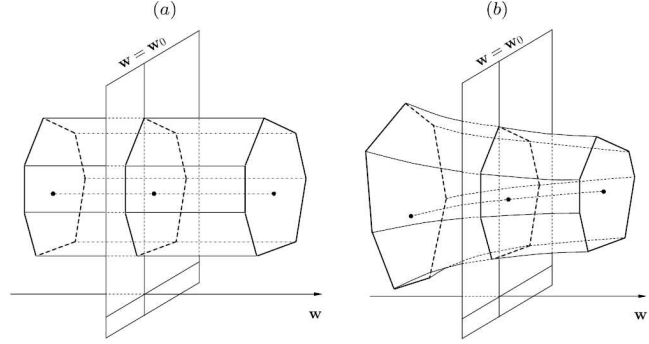}

\caption{Construction of (a) the local constant and (b) the local
bilinear $\tau
$-quantile regions as described in Sections \protect\ref{locconstsec}
and \protect\ref
{loclinsec}.
%from the intersection (over all possible directions $\bf u$) of the
%upper quantile halfspaces associated with the
%constrained-to-be-parallel-to-the-space-of-covariates (a) or
%unconstrained (b) $(\tau\mathbf{u})$-quantile hyperplanes.
}%
\label{locconstrucpic}
\end{figure}

%If we only care about %$\mathbf{w}_0$-conditional contours, that is,
% $\mathbf{w}_0$-cuts for a few selected number of $\mathbf{w}_0$
%values,
The above weighting scheme can be applied in the computation of the
weighted cylindrical regions generated by the hyperplanes in (\ref
{scalara})%(that are parallel to the space of covariates)
; more precisely, these cylindrical regions, with edges parallel to the
space of covariates, are obtained by computing the intersection (over
all $\mathbf{u}$'s, for fixed $\tau$) of the upper quantile
halfspaces associated with the quantile hyperplanes in (\ref
{scalara}); see Figure~\ref{locconstrucpic}(a).

The intersection with the $\mathbf{w}=\mathbf{w}_0$ hyperplane of
these cylindrical regions yields a local constant estimate, $\partial
\hat{R}^{(n)\mathrm{const}}_{\mathbf{w}_0}(\tau)$ say, of the
corresponding population $\mathbf{w}_0$-cut $\partial{R}_{\mathbf
{w}_0}(\tau)$; see\vadjust{\goodbreak} Section~\ref{asympsec} for asymptotic results. Of
course, the resulting local constant $\tau$-quantile/depth contours, namely
\[
\partial\hat\mathbf{ R}^{(n)\mathrm{const}}(\tau) := \bigcup
_{\mathbf{w}_0\in\R^{p-1}}\partial\hat {R}^{(n)\mathrm{const}}_{\mathbf{w}_0}(\tau),
\]
are not (globally) cylindrical, but rather adapt to the underlying
possibly nonlinear and/or heteroskedastic dependence structures.
%
%f4 #&#

This approach, which constitutes a generalization of the local constant
approach adopted elsewhere for single-output regression, has many
advantages. The main one is parsimony: each quantile hyperplane
involved in the construction %of the weighted contours
only entails $m$ parameters,
which is strictly less than in the local bilinear approach of
%%described in
the next section. On the other hand, the local constant approach
does not provide any information on, nor does take any advantage of,
the behavior of $\mathbf{w}$-cuts for~$\mathbf{w}$ values in
the neighborhood of $\mathbf{w}_0$, and its boundary performances are
likely to be poor. These two reasons, in traditional contexts, have
motivated the development of local linear and local polynomial methods;
see \cite{FG96} for a classical reference. Local linear methods were
successfully used in single-output quantile regression (\cite
{YJ97,YJ98,YL04,ZW09}). Considering them in the present context, thus,
is a quite natural idea.

%s4.3 #&#
\subsection{Local bilinear quantile/depth contours}
\label{loclinsec}

Assume that the distribution of $(\mathbf{W}^{\prime},\mathbf
{Y}^{\prime})^{\prime}$
is smooth enough that the coefficients of $\mathbf{w}$-conditional
quantile hyperplanes are differentiable with respect to $\mathbf{w}$.
%Denote by $(\mathbf{w}_0\pr, \mathbf{y})$ a point of the
%nonparametric contour $ \partial\mathbf{ R}(\tau)$.
Getting back to the first characterization (\ref{ddprems}) and (\ref
{unweightedfirst}) of quantile hyperplanes, the (restricted) $\mathbf
{w}_0$-conditional $\boldsymbol{\tau}$-quantile hyperplane of
$\mathbf{Y}$
defined in (\ref{condquantrestr}) and (\ref{CondOpt}) has equation (in
$\mathbf{y}$ -- of course, in $\mathbf{w}$, we just have $\mathbf
{w}=\mathbf{w}_0$)
%
%e4.2 #&#
\begin{eqnarray}
\label{linapprox.lc} \mathbf{u}^{\prime}\mathbf{y} - \bigl(a_{{\boldsymbol{\tau
}};\mathbf{w}_0} ,
\mathbf{c}^{\prime}_{{\boldsymbol{\tau}};\mathbf{w}_0} \bigr) \left( %
\begin{array} {c}1
\\
%\mathbf{w}_0 \\
{\boldsymbol{\Gamma}}^{\prime}_\mathbf{u}\mathbf{y}
\end{array} %
\right) = 0. %, \mbox{ with } \mathbf{a}^{(\mathbf{w}_0)\prime}_{\boldsymbol{\tau}}
%=: (a^{(
%\mathbf{w}_0)}_{{\boldsymbol{\tau}};0} , \mathbf{a}^{(\mathbf{w}_0)
%\prime}_{{
%\boldsymbol{\tau}};1}),
\end{eqnarray}
The same hyperplane equation, relative to a point $\mathbf{w}$ in the
neighborhood of $\mathbf{w}_0$, takes the form
%
%e4.3 #&#
\begin{eqnarray}
\label{linapprox} &&\mathbf{u}^{\prime}\mathbf{y} - \bigl(a_{{\boldsymbol{\tau}};\mathbf{w}_0} ,
\mathbf{c}^{\prime}_{{\boldsymbol{\tau}};\mathbf{w}_0} \bigr) \lleft( %
\begin{array} {c}1
\\
% \\ \mathbf{w}_0 \\
\boldsymbol{\Gamma}^{\prime}_\mathbf{u}\mathbf{y}
\end{array} %
\rright)
\nonumber
\\[-8pt]
\\[-8pt]
&&\quad {} -(\mathbf{w} - \mathbf{w}_0)^{\prime} \bigl(\dot
\mathbf {a}_{{\boldsymbol{\tau}
};\mathbf{w}_0} , %+ \mathbf{a}^{(\mathbf{w}_0) }_{{\boldsymbol{\tau}}
%;1}
%, \dot{\mathbf{a}}^{(\mathbf{w}_0) \prime}_{{\boldsymbol{\tau}}
%;1},
 \dot { \mathbf{c}}^{\prime}_{{\boldsymbol{\tau}};\mathbf{w}_0} \bigr) \lleft(
\begin{array} {c}1
\\
%\mathbf{w}_0 \\
{\boldsymbol{\Gamma}}^{\prime}_\mathbf{u}\mathbf{y}
\end{array} %
\rright) + \mathrm{o}\bigl(\Vert \mathbf{w} -
\mathbf{w}_0 \Vert\bigr)=0,
\nonumber
\end{eqnarray}
where $\dot\mathbf{a}_{{\boldsymbol{\tau}};\mathbf{w}_0}$ \label{defac}
stands for the gradient of $\mathbf{w}\mapsto{ a}_{{\boldsymbol{\tau}
};\mathbf{w}}$ and $\dot{ \mathbf{c}}_{{\boldsymbol{\tau}};\mathbf{w}}$
for the Jacobian matrix of $\mathbf{w}\mapsto{ \mathbf
{c}}_{{\boldsymbol{\tau}};\mathbf{w}}$, respectively, both taken at
$\mathbf{w}= \mathbf{w}_0$.
In order to express this equation into the equivalent quantile
formulation in (\ref{dd}) and (\ref{unweighted}), note that
we have $\mathbf{b}_{{\boldsymbol{\tau}};\mathbf{w}_0}=\mathbf{u}
- {\boldsymbol{\Gamma}}_\mathbf{u}\mathbf{c}_{{\boldsymbol{\tau
}};\mathbf{w}_0}$, which
entails $\dot{\mathbf{b}}_{{\boldsymbol{\tau}};\mathbf
{w}_0}=-{\boldsymbol{\Gamma}
}_\mathbf{u} \dot{ \mathbf{c}}_{{\boldsymbol{\tau}};\mathbf
{w}_0}$, where
$\dot{\mathbf{b}}_{{\boldsymbol{\tau}};\mathbf{w}_0}$ is the Jacobian
matrix of $\mathbf{w}\mapsto{ \mathbf{b}}_{{\boldsymbol{\tau
}};\mathbf
{w}}$ at $\mathbf{w}= \mathbf{w}_0$. Neglecting the $\mathrm{o}(\Vert\mathbf
{w} - \mathbf{w}_0 \Vert)$ term, (\ref{linapprox}) then rewrites,
after some algebra, as
%
%e4.4 #&#
\begin{eqnarray}
\label{linapprox'} && \bigl( {\mathbf{b}}^{\prime}_{{\boldsymbol{\tau}};\mathbf{w}_0} -
\mathbf{w}_0^{\prime} \dot{\mathbf{b}}^{\prime}_{{\boldsymbol{\tau
}};\mathbf{w}_0}
\bigr) \mathbf{y}
\nonumber
\\[-8pt]
\\[-8pt]
& &\quad {} - \bigl( a_{{\boldsymbol{\tau}};\mathbf{w}_0} - \mathbf{w}_0^{\prime
}
\dot \mathbf {a}_{{\boldsymbol{\tau}};\mathbf{w}_0 } , \dot{\mathbf{a}}^{\prime
}_{{\boldsymbol{\tau}};\mathbf{w}_0}
, - (\vec \dot{\mathbf {c}}_{{\boldsymbol{\tau}
};\mathbf{w}_0})^{\prime} \bigr) \lleft(
\begin{array} {c}1
\\
\mathbf{w}
\\
\mathbf{w} \otimes \bigl({\boldsymbol{\Gamma}}_\mathbf{u}^{\prime
}
\mathbf{y} \bigr) \end{array} %
\rright) =0.
\nonumber
\end{eqnarray}
Letting $\bar\mathbf{x}:=(1,\bar\mathbf{w}^{\prime})^{\prime
}:=(1,\mathbf
{w}^{\prime}, (\mathbf{w}\otimes{\boldsymbol{\Gamma}}_\mathbf
{u}^{\prime}\mathbf
{y})^{\prime})^{\prime}$, the latter equation is of the form
$
{\boldsymbol{\beta}}^{\prime}_{\boldsymbol{\tau}} \mathbf
{y}-{\boldsymbol{\alpha}}^{\prime}_{\boldsymbol{\tau}} (1,\bar
\mathbf{w}^{\prime})^{\prime}=0$,
with
${\boldsymbol{\beta}}^{\prime}_{\boldsymbol{\tau}} \mathbf{u} =
(
{\mathbf{b}}^{\prime}_{{\boldsymbol{\tau}};\mathbf{w}_0}
-
\mathbf{w}_0^{\prime}\dot{\mathbf{b}}^{\prime}_{{\boldsymbol{\tau
}};\mathbf{w}_0}
)
\mathbf{u} =
{\mathbf{b}}^{\prime}_{{\boldsymbol{\tau}};\mathbf{w}_0}
\mathbf{u} =1
$
since $ \dot{\mathbf{b}}_{{\boldsymbol{\tau}};\mathbf
{w}_0}^{\prime}\mathbf{u}
= \break -\dot{ \mathbf{c}}_{{\boldsymbol{\tau}};\mathbf{w}_0}^{\prime
}{\boldsymbol{\Gamma}
}_\mathbf{u}^{\prime}\mathbf{u}
=0$.
% -- provided, however, that the vector of regressors $\mathbf{x}=(1,
%\mathbf{w}\pr)\pr$ is augmented to $(1,\bar\mathbf{w}\pr)\pr=$.
Comparing\vspace{2pt} with (\ref{dd}), this suggests a local linear approach based
on weighted quantile hyperplanes (in the $mp$-dimensional
regressor-response space associated with the augmented regressor $\bar
\mathbf{x}$, that is, the $(\bar\mathbf{w}^{\prime},\mathbf
{y}^{\prime})^{\prime}
$-space), yielding weighted empirical quantile hyperplanes with equations
%
%e4.5 #&#
\begin{equation}
\label{augmentedempiric} {\boldsymbol{\beta}}^{(n)\prime}_{{\boldsymbol{\tau}};\omega
}\mathbf{y} -{
\boldsymbol{\alpha}}^{(n)\prime}_{{\boldsymbol{\tau}};\omega
} \bigl(1,\bar\mathbf
{w}^{\prime} \bigr)^{\prime}=0, % \mbox{ where } \bar\mathbf{x}\pr= (1, \mathbf{w}, \mathbf{w}\otimes
%\mathbf{y})
\end{equation}
based on the same sequences of weights $\omega^{(n)}_i:=\omega^{(n)}
_{\mathbf{w}_0}(\mathbf{W}_i)$, $i=1,\ldots, n$, as in Section~\ref{weightdefsec}. Interpretation of the results, however, is easier from
(\ref{linapprox}) than from (\ref{linapprox'}).
The left-hand side of (\ref{linapprox}) indeed splits naturally into
two parts of independent interest: (i) the first one, made of the first
two terms, yields the equation of the $\mathbf{w}_0$-conditional
$\boldsymbol{\tau}$-quantile hyperplane of $\mathbf{Y}$, hence
provides the
required information for constructing the empirical $\mathbf
{w}_0$-cuts, whereas (ii) the second part (the third term)
provides the linear (linear with respect to $(\mathbf{w}-\mathbf
{w}_0)$; actually, bilinear in $(\mathbf{w}-\mathbf{w}_0)$ and
${\boldsymbol{\Gamma}}_\mathbf{u}^{\prime}\mathbf{y}$) correction
required for a small
perturbation $(\mathbf{w}-\mathbf{w}_0)$ of the value of the
conditioning variable. Therefore, the important quantities to be
recovered from ${\boldsymbol{\alpha}}^{(n)}_{{\boldsymbol{\tau
}};\omega}$ and ${\boldsymbol{\beta}}^{(n)}_{{\boldsymbol{\tau
}};\omega}$ are estimations of these two
parts, which are easily obtained by
\begin{enumerate}[(ii)]
\item[(i)] letting $\mathbf{w}=\mathbf{w}_0$ in (\ref
{augmentedempiric}), which yields the equation
\[
{\boldsymbol{\beta}}^{(n)\prime}_{{\boldsymbol{\tau}};\omega
}\mathbf{y} -{\boldsymbol{
\alpha}}^{(n)\prime}_{{\boldsymbol{\tau}};\omega
} \bigl(1, \mathbf{w}_0^{\prime}
, \bigl(\mathbf{w}_0 \otimes{\boldsymbol{\Gamma}}_\mathbf{u}^{\prime
}
\mathbf{y} \bigr)^{\prime} \bigr)^{\prime}=0
\]
of an empirical hyperplane providing an estimate of the two first terms
in (\ref{linapprox}), namely, the $\mathbf{w}_0$-conditional
$\boldsymbol{\tau}$-quantile hyperplane;
\item[(ii)] subtracting the latter equation
from (\ref{augmentedempiric}), which provides the bilinear correction term.
\end{enumerate}
The bilinear nature of the local approximation in (ii) is easily
explained by the fact that, in general, unless the $\mathbf
{w}_0$-conditional and $\mathbf{w}$-conditional $\boldsymbol{\tau}$-quantile
hyperplanes are parallel to each other, no higher-dimensional
hyperplane can run through both (for instance, two mutually skew
non-intersecting straight lines in $\R^3$ do not span a plane).
Omitting the additional $\mathbf{W}\otimes({\boldsymbol{\Gamma
}}_\mathbf
{u}^{\prime}\mathbf{Y})$ regressors (in (i) above) may result in
inconsistent estimators of the $\mathbf{w}_0$-conditional $\boldsymbol
{\tau}
$-quantile hyperplanes. The resulting regions in $\R^{m+p-1}$, are not
polyhedral anymore, but delimited by ruled quadrics (hyperbolic
paraboloids for $m=2$ and $p-1=1$), the intersections of which with the
$\mathbf{w}=\mathbf{w}_0$ hyperplane yield polyhedral estimated
$\mathbf{w}_0$-cuts; see Figure~\ref{locconstrucpic}(b).

The local bilinear approach is more informative than the local constant
one, and should be more reliable at boundary points; the price to be
paid is an increase of the covariate space dimension (due to the
presence of the regressors $\mathbf{W}$ and $\mathbf{W}\otimes
({\boldsymbol{\Gamma}}_\mathbf{u}^{\prime}\mathbf{Y})$ in (\ref
{augmentedempiric})),
hence of the number of free parameters ($mp$ %free parameters
instead of $m$ for the local constant method). Note however that the
smoothing features of the problem, namely the dimension of kernels,
remains unaffected ($p-1$, irrespective of $m$).

%For asymptotic results and a discussion of bandwidth selection, see
%Section~3.4.

%s5 #&#
\section{Asymptotics}
\label{asympsec}
% and bandwidth choice}

Throughout this section, we fix $\mathbf{w}_0$ and $\boldsymbol{\tau
}=\tau
\mathbf{u}$, hence also $a_{\boldsymbol{\tau};\mathbf{w}_0}$ and
$\mathbf
{c}_{\boldsymbol{\tau};\mathbf{w}_0}$, and write, for %notational
simplicity,
%$\mathbf{Y}(\mathbf{u}):=({Y}_1(\mathbf{u}),\ldots,{Y}_m(\mathbf{u}))
%\pr:=({Y}_\mathbf{u}, \mathbf{Y}_{\mathbf{u}}^{\perp\prime})\pr$, with
${Y}_\mathbf{u}:=\mathbf{u}^{\prime}\mathbf{Y}$ and $\mathbf
{Y}_{\mathbf
{u}}^{\perp} :={\boldsymbol{\Gamma}}_\mathbf{u}^{\prime}\mathbf
{Y}$. %
Asymptotic results require some regularity assumptions on the density $f$,
the kernel $K$, and the bandwidth $h_n$.

\renewcommand{\theAssumption}{(A\arabic{Assumption})}
\begin{Assumption}\label{as1}
\begin{enumerate}[(iii)]
\item[(i)] The $n$-tuple $(\mathbf{W}_i^{\prime},
\mathbf
{Y}_i^{\prime})^{\prime}$, $i=1,\ldots, n$ is an \mbox{i.i.d.}
sample from
$(\mathbf{W}^{\prime}, \mathbf{Y}^{\prime})^{\prime}$.

\item[(ii)] The density $\mathbf{w}\mapsto f^\mathbf{W}(\mathbf{w})$ of
$\mathbf{W}$ is continuous and strictly positive at $\mathbf{w}_0$.

\item[(iii)] For any $\mathbf{t}\in\R^{m-1}$, there exist a neighborhood
$B_\mathbf{t}$ of $a_{\taub;\mathbf{w}_0} + \mathbf{c}_{\taub
;\mathbf{w}_0}^{\prime}\mathbf{t}$ and a neighborhood $\mathbf
{B}_\mathbf
{t}(\mathbf{w}_0)$ of $\mathbf{w}_0$ such that
$s\mapsto f^{{Y}_\mathbf{u}|\mathbf{Y}_{\mathbf{u}}^{\perp} =
\mathbf{t},\mathbf{W} = \mathbf{w}}(s)$ is continuous over $s\in
B_\mathbf{t}$, uniformly in $\mathbf{w}\in\mathbf{B}_\mathbf
{t}(\mathbf{w}_0)$, and $\mathbf{w}\mapsto f^{{Y}_\mathbf{u}|\mathbf
{Y}_{\mathbf{u}}^{\perp} = \mathbf{t},\mathbf{W} = \mathbf{w}}(s)$
is continuous over $\mathbf{w}\in\mathbf{B}_\mathbf{t}(\mathbf
{w}_0)$ for all $s\in B_\mathbf{t}$.

\item[(iv)] The density $f^{\mathbf{Y}_{\mathbf{u}}^{\perp}\vert\mathbf
{W} = \mathbf{w}}(\mathbf{t})$ of $\mathbf{Y}_{\mathbf{u}}^{\perp
}$ conditional on $\mathbf{W} = \mathbf{w}$ is continuous with
respect to $\mathbf{w}$ over a neighborhood of $\mathbf{w}_0$, except
perhaps for a set of $\mathbf{t}$ values of $f^{\mathbf{Y}_{\mathbf
{u}}^{\perp}}$-measure zero.

\item[(v)] The $m\times m$ matrix
\begin{eqnarray*}
%\label{fstar}
\mathbf{G}_{\taub;\mathbf{w}_0} := \int_{\R^{m-1}}
\lleft( %
\begin{array} {c@{\quad}c} 1 &
\mathbf{t}^{\prime}
\\
\mathbf{t} & \mathbf{t}\mathbf{t}^{\prime} \end{array} %
\rright) %
f^{{Y}_\mathbf{u} \vert\mathbf{Y}_{\mathbf{u}}^{\perp}=\mathbf
{t}, \mathbf{W}=\mathbf{w}_0} \bigl(a_{\taub;\mathbf{w}_0} + \mathbf
{c}_{\taub;\mathbf{w}_0}^{\prime}\mathbf{t} \bigr) f^{\mathbf{Y}_{\mathbf{u}}^{\perp} \vert\mathbf{W}=\mathbf
{w}_0}(\mathbf{t})
\, \mathrm{d}\mathbf{t}
\end{eqnarray*}
is finite and positive definite.
\end{enumerate}
\end{Assumption}
\begin{Assumption}\label{as2} The kernel function $K$
\begin{enumerate}[(iii)]
\item[(i)] is a compactly
supported bounded probability density over $\R^{p-1}$ %$S_K$, say),
such that
\item[(ii)] $\int_{\R^{p-1}}\mathbf{w}K(\mathbf{w})\, \mathrm{d}\mathbf{w}
=\mathbf{0}$ and ${\boldsymbol{\mu}}^K_2 :=\int_{\R^{p-1}}\mathbf
{ww}^{\prime}
K(\mathbf{w}) \,\mathrm{d}\mathbf{w}$ is  positive definite.
\end{enumerate}
\end{Assumption}
\begin{Assumption}\label{as3} The bandwidth $h_n$ is such that
$\lim_{\ny}h_n = 0$ and $\lim_{\ny}nh_n^{p-1} = \infty$.
\end{Assumption}

The conditions we are imposing in Assumption~\ref{as1} are quite mild. For
example, Assumption~\ref{as1}(ii) is the same as Condition
(A)(iii) in \cite{Faetal94} and Assumption (A1)(i) in \cite{Hetal09};
Assumption~\ref{as1}(iii)--(v) are similar to
Condition (A)(i, iv) in \cite{Faetal94} and Condition (A1)(ii) in
\cite{Hetal09}, where
the existence and positive-definiteness ensure the invertibility of
$\mathbf{G}_{\taub;\mathbf{w}_0}$ in Theorem~\ref{bahadur}.

Assumptions \ref{as2} and \ref{as3} on the kernel function
and the bandwidth also are quite standard in the nonparametric literature.
For example, any compactly supported symmetric density function
satisfies Assumption~\ref{as2}. The compact support of $K$ in Assumption~\ref{as2} is only a technical assumption to simplify the proof of theorems.
In practice, Gaussian kernels can be considered; indeed, at the cost of more
involved proof, the compact support assumption in Theorems \ref
{bahadur} and \ref{a.normality} can be replaced with the assumption that both
$C_0^K:=\int_{\mathbb{R}^{p-1}} K^2(\mathbf{w}) \,\mathrm{d}\mathbf{w}$
and
$C^K_2:=\int_{\mathbb{R}^{p-1}} \mathbf{w}\mathbf{w}' K^2(\mathbf
{w}) \,\mathrm{d}\mathbf{w}$ are finite.
As for Assumption~\ref{as3}, it is the usual one in the i.i.d.
setting; see Section~\ref{Bandsec} for a discussion.

%\subsubsection{Local constant case}
%\label{asymlocconstsec}

%We first give the asymptotic results for the local constant case.
Let
$\XX_{\mathbf{u}}^{c} := (1,\mathbf{Y}_{\mathbf{u}}^{\perp\prime}
)^{\prime}$
and
$\XX^\ell_\mathbf{u} := (1,\mathbf{Y}_{\mathbf{u}}^{\perp\prime
})^{\prime}\otimes(1, (\mathbf{W} - \mathbf{w}_0)^{\prime
})^{\prime}
$,
where the superscript $c$ and $\ell$ stand for the local constant and
local bilinear cases, respectively. For $(\mathbf{W}, \mathbf{Y}) =
(\mathbf{W}_i, \mathbf{Y}_i)$,
we use the notation ${Y}_{i\mathbf{u}}$, $\mathbf{Y}_{i\mathbf
{u}}^{\perp}$, $\XX_{i\mathbf{u}}^{c}$, $\XX_{i\mathbf{u}}^{\ell
}$, etc. in an obvious way.

Referring to (\ref{linapprox.lc}) for the notation, the parameter of
interest for the local constant case is $\thetab^c=\thetab
^c_{{\boldsymbol{\tau}};\mathbf{w}_0} :=(a_{{\boldsymbol{\tau
}};\mathbf{w}_0},\mathbf
{c}_{{\boldsymbol{\tau}};\mathbf{w}_0}^{\prime})^{\prime}$,
whereas, in the local
bilinear case (see (\ref{linapprox})), we rather have to estimate
%
%e5.1 #&#
\begin{equation}
\label{paramlin} \thetab^\ell =\thetab^\ell_{{\boldsymbol{\tau}};\mathbf{w}_0} :=
\vec \lleft( %
\begin{array} {c@{\quad}c} a_{{\boldsymbol{\tau}};\mathbf{w}_0} &
\mathbf{c}_{{\boldsymbol
{\tau}};\mathbf
{w}_0}^{\prime}
\\
\dot\mathbf{a}_{{\boldsymbol{\tau}};\mathbf{w}_0} & \dot\mathbf {c}_{{\boldsymbol{\tau}};\mathbf{w}_0}^{\prime}
\end{array} %
\rright).
\end{equation}
The local constant and local bilinear methods described in the previous
sections provide estimators of the form
$\hat{\thetab}^{c(n)} :=(\hat{a} ,\hat{\mathbf{c}}^{\prime
})^{\prime}$ and
%
%e5.2 #&#
\begin{equation}
\label{Zudi1.1a} \hat{\thetab}^{\ell(n)} := \vec \lleft( %
\begin{array} {c@{\quad}c} \hat{a} & \hat{\mathbf{c}}^{\prime}
\\
\hat{\dot\mathbf{a}} & \hat{\dot\mathbf{c}}^{\prime} \end{array}
\rright)
\end{equation}
(we should actually discriminate between $(\hat{a}, \hat\mathbf
{c}')=(\hat{a}^c, \hat\mathbf{c}^{c\prime})$ and $(\hat{a}, \hat
\mathbf{c}')=(\hat{a}^\ell, \hat\mathbf{c}^{\ell\prime})$, but
will not do so in order to avoid making the notation too heavy); those
estimators are defined as the corresponding minimizer $\thetab^r$ of
%
%e5.3 #&#
\begin{equation}
\label{Zudi1.1.lc} \sum_{i=1}^n
K_h(\mathbf{W}_i - \mathbf{w}_0)
\rho_\tau \bigl(Y_{i\mathbf{u}} - \thetab^{r\prime}
\XX_{i\mathbf{u}}^{r} \bigr) ,\qquad r=c,\ell .
\end{equation}

The following result provides Bahadur representations for $\hat
{\thetab}
^{c(n)}$ and $\hat{\thetab}^{\ell(n)}$.

%th5.1 #&#
\begin{Theor}[(Bahadur representations)]
\label{bahadur}
%\label{bahadur.lc}
Let Assumptions
\textup{\ref{as1}}, \textup{\ref{as2}(i)} and \textup{\ref{as3}} hold, assume that $\mathbf{w}\mapsto
(a_{{\boldsymbol{\tau}};\mathbf{w}} ,\mathbf{c}_{{\boldsymbol{\tau
}};\mathbf
{w}}^{\prime})^{\prime}$ is continuously differentiable at
$\mathbf{w}_0$,
and write $\psi_{\tau}(y) :={\tau}-I[ y < 0] $.
%, with gradient $\left(\dot{a}^{(\mathbf{x}_0)}_{{\boldsymbol{\tau}}} ,
%\dot{
%\mathbf{c}}^{(\mathbf{x}_0)\prime}_{\boldsymbol{\tau}}\right)$.
Then, as $\ny$,
\begin{eqnarray}
\label{bahadur.1} && \sqrt{n h^{p-1}_n}
\mathbf{M}_h^r \bigl( \hat{\boldsymbol{
\theta}}^{r(n)} -{\boldsymbol{\theta }}^r \bigr)
\nonumber
\\[-8pt]
\\[-8pt]
&&\quad = \frac{
{\boldsymbol{\eta}}^r_{\taub;\mathbf{w}_0}}{\sqrt{n h_n ^{p-1}}} \sum_{i=1}^n
K \biggl(\frac{\mathbf{W}_{i} -\mathbf{w}_0}{h_{n }} \biggr) \psi_{\tau} \bigl(
{Z}_{i\mathbf{u}}^r({ \boldsymbol{\theta}}) \bigr) \bigl(
\mathbf{M}_h^r \bigr)^{-1}
\XX^r_{i\mathbf{u}} +\mathrm{o}_\mathrm{P}(1) ,
\nonumber
%\label{bahadur.1.lc}
\end{eqnarray}
where
${Z}_{i\mathbf{u}}^r({\thetabv})
:={Y}_{i\mathbf{u}} -{\thetabv}^{\prime}\XX_{i\mathbf{u}}^{r}
$ $(r=c,\ell)$,
$\mathbf{M}_h^c :=\mathbf{I}_m$,
$\mathbf{M}_h^\ell:=\mathbf{I}_m\otimes\diag(1, h_n\mathbf{I}_{p-1})$,
\begin{eqnarray*}
%\label{etatau}
%\label{etatau.lc}
{\boldsymbol{\eta}}^{c}_{{\boldsymbol{\tau}};\mathbf{w}_0}:=
%\frac{1}
\bigl({f^\mathbf{W}(\mathbf{w}_0)}
\bigr)^{-1} \mathbf{G}_{\taub;\mathbf{w}_0}^{-1} \quad \mbox{and}\quad  {
\boldsymbol{\eta}}^\ell_{\taub;\mathbf{w}_0} := %\frac{1}
\bigl({f^\mathbf{W}(\mathbf{w}_0)} \bigr)^{-1}
\mathbf{G}_{\taub;\mathbf{w}_0}^{-1} \otimes \diag \bigl(1, \bigl({
\boldsymbol{\mu}}^K_2 \bigr)^{-1} \bigr),
\end{eqnarray*}
with $\mathbf{G}_{\taub;\mathbf{w}_0}$ defined in Assumption \textup{\ref{as1}(v)}
(the result for the local constant case does not require \textup{\ref{as2}(ii)}).
\end{Theor}

This result, along with Assumption~\ref{as4}
below, entails the asymptotic normality of $\hat{\thetab}^{r(n)}$,
$r=c,\ell$. That assumption deals with the existence, at $\mathbf{w}
= \mathbf{w}_0$, of the second derivatives of $\mathbf{w}\mapsto
(a_{{\boldsymbol{\tau}};\mathbf{w}} ,\mathbf{c}_{{\boldsymbol{\tau
}};\mathbf
{w}}^{\prime})^{\prime}$.
With $\mathbf{c}_{{\boldsymbol{\tau}};\mathbf{w}}
=:({c}_{{\boldsymbol{\tau}
};\mathbf{w},1}, \ldots, {c}_{{\boldsymbol{\tau}};\mathbf
{w},m-1})^{\prime}$,
denote by $\dot\mathbf{a}_{{\boldsymbol{\tau}};\mathbf{w}}$
and $\dot\mathbf{c}_{{\boldsymbol{\tau}};\mathbf{w},j}$ the
$(p-1)\times1$
vectors of first derivatives and by $\ddot\mathbf{a}_{{\boldsymbol
{\tau}
};\mathbf{w}}$ and $\ddot\mathbf{c}_{{\boldsymbol{\tau}};\mathbf{w},j}$
the $(p-1)\times(p-1)$ matrices of second derivatives (when they
exist) of $\mathbf{w}\mapsto a_{{\boldsymbol{\tau}};\mathbf{w}}$ and
$\mathbf{w}\mapsto
{c}_{{\boldsymbol{\tau}};\mathbf{w},j}$, respectively
(recall that $\dot\mathbf{a}_{{\boldsymbol{\tau}};\mathbf{w}}$ and
$\dot\mathbf{c}_{{\boldsymbol{\tau}};\mathbf{w}}=
(\dot\mathbf{c}_{{\boldsymbol{\tau}};\mathbf{w},1},
\ldots,
\dot\mathbf{c}_{{\boldsymbol{\tau}};\mathbf{w},m-1})^{\prime}$
%(the Jacobian matrix of $\mathbf{w}\mapsto{c}_{{\boldsymbol{\tau}};
%\mathbf{w}}$)
were already defined in page \pageref{defac}).
Finally, write $\ddot\mathbf{c}_{{\boldsymbol{\tau}};\mathbf
{w}}^{\prime}$ for
the $(p-1)\times(m-1)(p-1)$ matrix
$(\ddot\mathbf{c}_{{\boldsymbol{\tau}};\mathbf{w},1},\ldots,\ddot
\mathbf
{c}_{{\boldsymbol{\tau}};\mathbf{w},m-1})$.

\begin{Assumption}\label{as4}
\begin{enumerate}[(iii)]
\item[(i)] The function
$\mathbf{w}\mapsto(a_{{\boldsymbol{\tau}};\mathbf{w}} ,\mathbf
{c}_{{\boldsymbol{\tau}};\mathbf{w}}^{\prime})^{\prime}$ is twice
continuously
differentiable at $\mathbf{w} = \mathbf{w}_0$, that is, $\ddot
\mathbf{a}_{{\boldsymbol{\tau}};\mathbf{w}}$ and $\ddot\mathbf
{c}_{{\boldsymbol{\tau}};\mathbf{w}}$ exist
in a neighborhood of $\mathbf{w}_0$
and are continuous with respect to $\mathbf{w}$ at $\mathbf{w}_0$.
\item[(ii)] The function
$\mathbf{w}\mapsto f^\mathbf{W}(\mathbf{w})$ is continuously
differentiable at $\mathbf{w} = \mathbf{w}_0$, that is, the
$(p-1)\times1$ vector of first derivatives of $f^\mathbf{W}$, $\dot
{f}^\mathbf{W}(\mathbf{w}) $, exists
in a neighborhood of $\mathbf{w}_0$
and is continuous with respect to $\mathbf{w}$ at $\mathbf{w}_0$.
\end{enumerate}
\end{Assumption}

The following matrices are involved in the
asymptotic bias and variance expressions of the asymptotic normality
result in Theorem~\ref{a.normality} below. Define
%
%e5.5 #&#
%e5.6 #&#
\begin{eqnarray}
\label{Sig(x).lc}\boldsymbol{\Sigma}_{\mathbf{w}}^{c} & := & \tau(1-
\tau) f^\mathbf{W}(\mathbf{w}) C^K_0 {\boldsymbol{
\eta}}_{{\boldsymbol{\tau}};\mathbf{w}}^{c} \biggl[ \int_{\mathbb{R}^{m-1}}
f^{\mathbf{Y}_\mathbf{u}^{\perp}\vert\mathbf{W} = \mathbf
{w}}(\mathbf{t}) %
\lleft( %
\begin{array}
{c@{\quad}c} 1& \mathbf{t}^{\prime}
\\
\mathbf{t} & \mathbf{t}\mathbf{t}^{\prime} \end{array} %
\rright) %
\,\mathrm{d}\mathbf{t} %\diag\Big(
%,\int_{\mathbb{R}^{p-1}}\mathbf{v}\mathbf{v}\pr K^2(\mathbf{v})d
%\mathbf{v}\Big)
\biggr] \boldsymbol{\eta}_{\boldsymbol{\tau};\mathbf{w}}^{c} ,
\\
%\end{equation}
%\begin{equation}
%\begin{eqnarray}
\label{Sig(x)}\boldsymbol{
\Sigma}^\ell_{\mathbf{w}} & := & \tau(1-\tau)f^\mathbf{W}(
\mathbf{w}) {\boldsymbol{ \eta}}^\ell_{{\boldsymbol{\tau}};\mathbf{w}}
\nonumber
\\[-8pt]
\\[-8pt]
& & {} \times \biggl[ \int_{\mathbb{R}^{m-1}} f^{\mathbf{Y}_\mathbf{u}^{\perp}\vert\mathbf{W} = \mathbf
{w}}(
\mathbf{t}) %
\lleft( %
\begin{array} {c@{\quad}c} 1&
\mathbf{t}^{\prime}
\\
\mathbf{t} & \mathbf{t}\mathbf{t}^{\prime} \end{array} %
\rright) %
\,\mathrm{d}\mathbf{t} \otimes \diag \bigl(
C^K_0,C^K_2 \bigr) \biggr] {
\boldsymbol{\eta}}^\ell_{{\boldsymbol{\tau}};\mathbf{w}} ,
\nonumber
\end{eqnarray}
%
%\end{equation}
and, for $r=c,\ell$,
%\begin{eqnarray}
%\mathbf{B}^{r}_\mathbf{w}
%:=
%f^\mathbf{W}(\mathbf{w})
%{\boldsymbol{\eta}}_{{\boldsymbol{\tau}};\mathbf{w}}^{r}
% \int_{\mathbb{R}^{m-1}}
%f^{{Y}_\mathbf{u}\vert\mathbf{Y}_\mathbf{u}^{\perp}=\mathbf{t},
%\mathbf{W} = \mathbf{w}} (a_{{\boldsymbol{\tau}};\mathbf{w}} +
%\mathbf{c}_{{\boldsymbol{\tau}};\mathbf{w}}\pr\mathbf{t})
% \nonumber
%\\
%& &
%
% \times f^{\mathbf{Y}_\mathbf{u}^{\perp}\vert\mathbf{W} = \mathbf{w}}(
%\mathbf{t} )
% \Big( \begin{array}{c} 1 \\
%\mathbf{t}
%\end{array}
% \Big)\otimes\Big[
%\mathbf{B}_{\mathbf{w};0}^{r}\Big( \begin{array}{c} 1 \\
% \mathbf{t}
%\end{array}
% \Big)
%\Black\Big] \mathrm{d}\mathbf{t},
%\label{B(x).lc}
%\end{eqnarray}
%%%%%%%
%
%e5.7 #&#
\begin{eqnarray}
\label{B(x).lc} \mathbf{B}^{r}_\mathbf{w} &:=
&f^\mathbf{W}(\mathbf{w}) { \boldsymbol{\eta}}_{{\boldsymbol{\tau}};\mathbf{w}}^{r}
\nonumber
\\[-8pt]
\\[-8pt]
& & {} \times \int_{\mathbb{R}^{m-1}} f^{{Y}_\mathbf{u}\vert\mathbf{Y}_\mathbf{u}^{\perp}=\mathbf{t},
\mathbf{W} = \mathbf{w}}
\bigl(a_{{\boldsymbol{\tau}};\mathbf{w}} +\mathbf {c}_{{\boldsymbol{\tau}};\mathbf{w}}^{\prime}\mathbf{t}
\bigr) f^{\mathbf{Y}_\mathbf{u}^{\perp}\vert\mathbf{W} = \mathbf
{w}}(\mathbf{t} ) \lleft( %
\begin{array} {c} 1
\\
\mathbf{t} \end{array} %
\rright)\otimes \biggl[ \mathbf{B}_{\mathbf{w};0}^{r}
\lleft( %
\begin{array} {c} 1
\\
\mathbf{t} \end{array} %
\rright) \biggr] \,\mathrm{d}\mathbf{t},
\nonumber
\end{eqnarray}
%
%%%%%%%%%
where
(putting $\ddot\mathbf{c}_{{\boldsymbol{\tau}};\mathbf
{w},0}:=\ddot\mathbf
{a}_{{\boldsymbol{\tau}};\mathbf{w}}$)
$\mathbf{B}_{\mathbf{w};0}^{c}$ is the $1\times m$ matrix
with $j$th entry
\[
B_{\mathbf{w};0,j}^{c}:= \tr \biggl[ \biggl( \ddot
\mathbf{c}_{{\boldsymbol{\tau}};\mathbf{w},j-1} +2 \frac{
\dot
\mathbf{c}_{{\boldsymbol{\tau}};\mathbf{w},j-1} (\dot{f}^{\mathbf
{W}}(\mathbf{w}))^{\prime}}{{f}^{\mathbf{W}}(\mathbf{w})} \biggr) {\boldsymbol{
\mu}}^K_2 \biggr], \qquad j=1,\ldots,m,
\]
and
$\mathbf{B}^\ell_{\mathbf{w};0}$ denotes the $p\times m$ matrix
with $(i,j)$th entry
\[
B^\ell_{\mathbf{w};0,ij}:= \tr \biggl[ \ddot\mathbf{c}_{{\boldsymbol{\tau}};\mathbf{w},j-1}
\int_{\mathbb{R}^{p-1}} w_{i-1} \mathbf{w} \mathbf{w}^{\prime}K(
\mathbf{w}) \,\mathrm{d}\mathbf{w} \biggr] , \qquad i=1,\ldots,p, j=1,\ldots,m ;
\]
here, we wrote $\mathbf{w}=(w_1,w_2,\ldots,w_{p-1})^{\prime}$,
$w_0=1$. We
then have:

%th5.2 #&#
\begin{Theor}[(Asymptotic normality)]
\label{a.normality}
%\label{a.normality.lc}
Let Assumptions \textup{\ref{as1}}--\textup{\ref{as4}} hold. Then, for $r=c,\ell$,
%
%e5.8 #&#
\begin{equation}
\sqrt{nh_n^{p-1}} \mathbf{M}_h^r
\biggl( \hat{\boldsymbol{\theta}}^{r(n)}-{\boldsymbol{
\theta}}^r-\frac
{h^2}{2} \mathbf {B}^r_{\mathbf{w}_0}
\biggr) \stackrel{\mathcal{L}} {\rightarrow} \mathcal{N} \bigl(\mathbf{0}, {
\boldsymbol{\Sigma}}^r_{\mathbf{w}_0} \bigr), \label{3.8a}
\end{equation}
as $n\to\infty$, where
$\stackrel{\mathcal{L}}{\rightarrow}$ denotes convergence in
distribution
(the result for the local bilinear case does not require Assumption \textup{\ref{as4}(ii)}).
\end{Theor}

%%%%%%%%%%%%%%%%%%%%%%%%%%%%%%%%%
\begin{Remark}
The local bilinear fitting has an expression of bias that is
independent of $\dot f^\mathbf{W}$. In contrast, the
local constant fitting has a large bias at the regions where the
derivative of
$f^\mathbf{W}$ is large, that is, it does not adapt to highly-skewed
designs (see \cite{FG96,FY03}).
%Comparing the local constant ($r=c$) and local linear ($r=\ell$)
%versions of (\ref{B(x).lc}) shows an obvious advantage of local linear
%fitting over the local constant approach is the reduction of bias
Another important advantage of local
bilinear fitting over the local constant approach is its much better
boundary behavior. This advantage often has been emphasized in the
usual regression settings when the regressors take values on a compact
subset of $\mathbb{R}^{p-1}$. For example, considering a univariate
random regressor $W$ ($p=2$) with bounded support ($[0, 1]$, say), it
can be proved, using an argument similar to the one developed in the
corresponding proof in \cite{FG96}, that asymptotic normality (with
the same rate) still holds at boundary points of the form $ch_{n}$,
where $c\in\R^+_0$, with asymptotic bias and variances of the same
form as in the local bilinear ($r=\ell$) versions of (\ref{B(x).lc})
and (\ref{Sig(x)}), with $p=2$, $\mathbf{w}_0$ replaced by ${w}_0=0^+$,
and $\int_{\mathbb{R}^{p-1}}$ by $\int_{-c}^{\infty}$;
see, for example, page 666 of \cite{Hetal09}.
\end{Remark}

\begin{Remark}
In practice, we may be concerned with the
estimation of the quantile regression functions at different
$\taub$'s simultaneously. Restricting to the estimation of $(\thetab
_{\taub_1;\mathbf{w}_0}^{\prime},\thetab_{\taub_2;\mathbf
{w}_0}^{\prime})^{\prime}
$, it can be shown by proceeding as in the proof of Theorem~\ref
{a.normality} that
$(\hat{\thetab}_{\taub_1;\mathbf{w}_0}^{\prime},\hat{\thetab
}_{\taub
_2;\mathbf{w}_0}^{\prime})^{\prime}$ is asymptotically normal with a
block-diagonal asymptotic covariance matrix, that is,
$\hat{\thetab}_{\taub_1;\mathbf{w}_0}$ and $\hat{\thetab}_{\taub
_2;\mathbf{w}_0}$
are asymptotically independent for $\taub_1\ne\taub_2$.
%, the argument for which is similar to that of Lemma~\ref{LemA6}.
\end{Remark}

%s6 #&#
\section{Bandwidth selection}\label{Bandsec}
While the choice of a kernel, as usual, has little impact on the final
result, selecting the bandwidth $h$ is more delicate. A full plug-in
estimator in principle could be derived from the asymptotic normality
result of Theorem~\ref{a.normality}, along the same lines as, for
instance, in Zhang and Lee \cite{ZL00}, who do it for mean
regression. Such an approach, however, requires the estimation of
several conditional densities, hence raises further problems, besides
being computationally quite heavy, certainly when several values of
$\tau$ are to be considered. A simpler heuristic rule is thus
preferable; the one we are describing here is adapted from \cite
{YL04}, where it is proposed in the context of single-output quantile
regression.

Without loss of generality, we restrict to
$p-1=1$ for notational simplicity, writing ${W}$ and ${w}$ for $\mathbf
{W}$ and $\mathbf{w}$, $h$ for $h_n$ and $\hat{\boldsymbol{\theta
}}_h=(\hat
{a}_{{\boldsymbol{\tau}};{ w}_0}^h, \hat{\mathbf{c}}_{{\boldsymbol
{\tau}};{
w}_0}^{h\prime})^{\prime}$ for the estimator of
${\boldsymbol{\theta}}=(a_{{\boldsymbol{\tau}};{w}_0}, \mathbf
{c}_{{\boldsymbol{\tau}
};{w}_0}^{\prime})^{\prime}$ associated with bandwidth $h$, respectively.
Throughout, the kernel $K$ is some symmetric density function, such
as the standard normal one. The objective is to minimize, with
respect to $h$, %\in[h_L, h_U]$ (with $0<h_L<h_U$ appropriately
%specified),
the asymptotic mean square error which, in view of Theorem~\ref{a.normality} with
$p-1=1$, after some straightforward algebra takes the form
%
%e6.1 #&#
\begin{equation}
\label{mseh} \mathit{MSE}(h)=\mathrm{E}(\hat{{\boldsymbol{\theta }}}_h-{{
\boldsymbol{ \theta}}})^{\prime} (\hat{{\boldsymbol{\theta}}}_h-{{
\boldsymbol{\theta}}}) \approx \frac{1} 4 h^4
B_{\boldsymbol{\tau}}^2+\frac{1} {nh} V_{\boldsymbol{\tau}},
\end{equation}
with
\[
B_{\boldsymbol{\tau}}^2:= \bigl(\mu_2^K
\bigr)^2 \Biggl(\ddot{a}_{{\boldsymbol
{\tau}};
w_0}^2+\sum
_{j=1}^{m-1} \ddot\mathbf{c}_{{\boldsymbol{\tau}};
w_0,j}^2
\Biggr) \quad \mbox{and} \quad %$$
%$$
V_{\boldsymbol{\tau}}:=
\frac{\tau(1-\tau)C_0^K}{ f^{W}({w}_0)} \tr \bigl(\mathbf{G}_{{\boldsymbol{\tau}};w_0}^{-1} \mathbf
{G}_{w_0} \mathbf{G}_{{\boldsymbol{\tau}};w_0}^{-1} \bigr),
\]
where $\ddot\mathbf{c}_{{\boldsymbol{\tau}}; w_0,j}$ is the second-order
derivative with respect to $w_0$ of the $j$th component of $\mathbf
{c}_{{\boldsymbol{\tau}}; w_0}$, $\mathbf{G}_{{\boldsymbol{\tau
}};w_0}$ is defined
in Assumption~\ref{as1}(v), and
\[
\mathbf{G}_{w_0}:=\int_{\mathbb{R}^{m-1}} f^{\mathbf{Y}_\mathbf{u}^{\perp}\vert{ W} = { w}_0}(
\mathbf{t}) %
\lleft( %
\begin{array} {c@{\quad}c} 1&
\mathbf{t}^{\prime}
\\
\mathbf{t} & \mathbf{t}\mathbf{t}^{\prime} \end{array} %
\rright) %
\,\mathrm{d}\mathbf{t}.
\]
The minimizer $h_{\boldsymbol{\tau}}$ of
%asymptotically optimal (interior) value,
%$h_{\tau}$, for $h$, following from
(\ref{mseh}) satisfies
%
%e6.2 #&#
\begin{equation}
\label{htau5} h_{\boldsymbol{\tau}}^{5}=\frac{V_{\boldsymbol{\tau}}}{n
B_{\boldsymbol{\tau}}^2}=
\frac
{\tau(1-\tau)C_0^K
\tr [\mathbf{G}_{{\boldsymbol{\tau}};w_0}^{-1} \mathbf
{G}_{w_0}
\mathbf{G}_{{\boldsymbol{\tau}};w_0}^{-1}  ]}{n (\mu
_2^K)^2f^{W}({w}_0)
(\ddot{a}_{{\boldsymbol{\tau}}; w_0}^2+\sum_{j=1}^{m-1} \ddot
\mathbf
{c}_{{\boldsymbol{\tau}},
w_0,j}^2 )},
\end{equation}
so that for any ${\boldsymbol{\tau}}_1, {\boldsymbol{\tau}}_2$,
%
%e6.3 #&#
\begin{equation}
\label{hratio} \biggl(\frac{h_{{\boldsymbol{\tau}}_1}}{h_{{\boldsymbol{\tau
}}_2}} \biggr)^{5} =
\frac{\tau_1
(1-\tau_1)}{\tau_2(1-\tau_2)} \frac{ (\ddot{a}_{{\boldsymbol
{\tau}}_2,
w_0}^2 +\sum_{j=1}^{m-1} \ddot\mathbf{c}_{{\boldsymbol{\tau}}_2,
w_0,j}^2 ) \tr (\mathbf{G}_{{\boldsymbol{\tau
}}_1,w_0}^{-1}
\mathbf{G}_{w_0}
\mathbf{G}_{{\boldsymbol{\tau}}_1,w_0}^{-1}  )}{ (\ddot
{a}_{{\boldsymbol{\tau}}_1,
w_0}^2 +\sum_{j=1}^{m-1} \ddot\mathbf{c}_{{\boldsymbol{\tau}}_1,
w_0,j}^2 )\tr (\mathbf{G}_{{\boldsymbol{\tau}
}_2,w_0}^{-1} \mathbf{G}_{w_0}
\mathbf{G}_{{\boldsymbol{\tau}}_2,w_0}^{-1}  )}.
\end{equation}

As in \cite{YL04}, we assume that $\ddot{a}_{\tau\mathbf{u}, w_0}$
and $\ddot\mathbf{c}_{\tau\mathbf{u}, w_0}$ do not depend on $\tau
$ (an assumption we do not make on ${a}_{\tau\mathbf{u},w_0}$ and
$\mathbf{c}_{\tau\mathbf{u},w_0}$).
If $ f^{{Y}_\mathbf{u}\vert\mathbf{Y}_\mathbf{u}^{\perp}=\mathbf
{t}, { W} = {w}_0} $ were a normal
density with mean $\mu_{\mathbf{t},w_0}$ and
variance $\sigma_{\mathbf{t},w_0}^{2}$, denoting by $\phi$ and $\Phi
$ the standard normal density and
distribution functions, respectively,
%and by
%$z_\tau=\Phi^{-1}(\tau)$ the standard normal quantile of order $\tau$,
we would have $f^{{Y}_\mathbf{u}\vert\mathbf{Y}_\mathbf{u}^{\perp
}=\mathbf{t}, {W} = {w}_0} (a_{{\boldsymbol{\tau}};{w}_0} +\mathbf
{c}_{{\boldsymbol{\tau}};{w}_0}^{\prime}\mathbf{t})=\sigma
_{\mathbf{t},{w}_0}^{-1}\phi
(\Phi^{-1}(\tau)) $, hence
\[
\mathbf{G}_{{\boldsymbol{\tau}};w_0}=\phi \bigl(\Phi^{-1}(1/2) \bigr)\int
_{\mathbb{R}^{m-1}} \sigma_{\mathbf{t},{w}_0}^{-1} f^{\mathbf{Y}_\mathbf{u}^{\perp}\vert{ W} = { w}_0}(
\mathbf{t}) %
\lleft( %
\begin{array} {c@{\quad}c} 1&
\mathbf{t}^{\prime}
\\
\mathbf{t} & \mathbf{t}\mathbf{t}^{\prime} \end{array} %
\rright) %
\,\mathrm{d}\mathbf{t}
\]
and
\[
\frac{\tr (\mathbf{G}_{{\boldsymbol{\tau
}}_1,w_0}^{-1} \mathbf
{G}_{w_0}
\mathbf{G}_{{\boldsymbol{\tau}}_1,w_0}^{-1} )}{\tr
(\mathbf{G}_{{\boldsymbol{\tau}}_2,w_0}^{-1}
\mathbf{G}_{w_0}
\mathbf{G}_{{\boldsymbol{\tau}}_2,w_0}^{-1} )}= \biggl[\frac{\phi(\Phi^{-1}(\tau_2))}{\phi(\Phi^{-1}(\tau
_1))} \biggr]^2.
\]
If we further assume that $\sigma_{\mathbf{t},w_0}^{2}=\sigma
_{w_0}^{2}$, (\ref{htau5}) for $\tau=1/2$ takes the form
%
%e6.4 #&#
\begin{equation}
\label{h1/2} h_{\mathbf{u}/2}^{5}=\frac{\uppi}{2} \biggl(
\frac{C_0^K}{n(\mu
_2^K)^2} \frac{\tr (\mathbf{G}_{w_0}^{-1} )
\sigma_{w_0}^2}{
f^W(w_0)  (\ddot{a}_{\mathbf{{u}}/2; w_0}^2+\sum_{j=1}^{m-1}
\ddot\mathbf{c}_{\mathbf{{u}}/2; w_0,j}^2 )} \biggr),
\end{equation}
while (\ref{hratio}) yields
%
%e6.5 #&#
\begin{equation}
\label{h1/2'} \biggl(\frac{h_{{\boldsymbol{\tau}}_1}}{h_{{\boldsymbol{\tau
}}_2}} \biggr)^{5} =
\frac{\tau_1(1-\tau_1)}{\tau_2(1-\tau_2)} \frac{(\phi(\Phi^{-1}(\tau_2)))^{2}}{(\phi(\Phi^{-1}(\tau_1)))^{2}}
\end{equation}
%
%which provides a practical way of modifying $h$ with $\tau$.
hence, for \({\boldsymbol{\tau}}_2={\mathbf{u}/2} \),
% \begin{equation}\label{htau}
$h_{{\boldsymbol{\tau}} }^{5}=(2/\uppi)\tau(1-\tau)
(\phi(\Phi^{-1}(\tau)))^{-2}h_{{\mathbf{u}/2}}^{5}
$.%\end{equation}

This latter expression still is not readily implementable. However,
(\ref{h1/2}) bears a strong relation to
%The remaining step consists in selecting a bandwidth $h_{{
%\mathbf{u}/2}}$ for $\tau=1/2$, that is, for median quantile (Least
%Absolute Deviation) hyperplanes. As in Yu and Lu (2004), we propose to
%adopt
% under some simple assumption. In
%fact,
% the automatic bandwidth $h_{1/2}$ can be expressed in terms of
the optimal bandwidth value $h_{\FZ}$ obtained by Fan and Zhang in
Theorem~1 of \cite{FZ08} for the estimation of the
conditional mean
% $\mathrm{E}\big[Y_\mathbf{u}-\mathbf{c}\pr\mathbf{Y}_\mathbf{u}^{
%\bot}-a({w_0}) |W=w_0\big]$ in the %. That optimal bandwidth is
% which is considered as follows:
%
% Denote $\epsilon_\mathbf{u}=Y_\mathbf{u}-a(W) - \mathbf{c}(W)\pr Y_
%\mathbf{u}^{\bot}$. We consider a
in the varying-coefficient linear regression
model
$Y_\mathbf{u}=a(W) + \mathbf{c}(W)^{\prime}\mathbf{Y}_\mathbf
{u}^{\bot
}+\epsilon_\mathbf{u}$ with $\operatorname{Var}(\epsilon_\mathbf{u} \mid W=w_0)=
\sigma^2_{w_0}$,
%=Var\{Y_\mathbf{u}|W=w_0 \}$,
namely % with $\sigma^2_{w_0}=Var\{Y_\mathbf{.
%Here as an approximation, we assume the conditional variance
%$\sigma_{\mathbf{t},w_0}^{2}$ of $\epsilon_\mathbf{u}$ given $(W=w_0,
%Y_\mathbf{u}^{\bot}=\mathbf{t})$ could be approximated by
%$\sigma_{w_0}^{2}$ that is independent of $\mathbf{t}$, as assumed in
%the usual varying-coefficient model (c.f., Fan and Zhang (1999)).
% It follows from Theorem~1 of Fan and Zhang (2008) that under
% local linear kernel fitting, the estimator for $(a({w_0}),
%\mathbf{c}({w_0})\pr)\pr$
% has an asymptotic optimal bandwidth (under MSE of (\ref{mseh})) given
%by
\[
h_{\FZ}^5=%\biggl(
\frac{C_0^K}{n(\mu_2^K)^2}
\frac{\tr (\mathbf
{G}_{w_0}^{-1} )\sigma_{w_0}^2}{
f^W(w_0)  (\ddot{a}_{w_0}^2+\sum_{j=1}^{m-1} \ddot\mathbf
{c}_{w_0,j}^2 )} = (2/\uppi)h_{\mathbf{{u}}/2}^5%\biggr)
.
\]
%
% \begin{equation}
%h_{\tau}^{5}=\frac{C_0^K \sigma^2_{w_0}
% \mathrm{tr}\left[\mathbf{G}_{{\boldsymbol{\tau}};w_0}^{-1}
%\mathbf{G}_{w_0}
% \mathbf{G}_{{\boldsymbol{\tau}};w_0}^{-1} \right]}{n (
%\mu_2^K)^2f^{W}({w}_0)
% \left(\ddot{a}_{{\boldsymbol{\tau}}; w_0}^2+\sum_{j=1}^{m-1} \ddot
%\mathbf{c}_{{
%\boldsymbol{\tau}}; w_0,j}^2\right)}
%\end{equation}
We therefore propose, for ${\boldsymbol{\tau}}=\tau\mathbf{{u}}$, the
bandwidth $ h_{\boldsymbol{\tau}}$ provided by
%
%e6.6 #&#
\begin{equation}
h_{\boldsymbol{\tau}}^5=\tau(1-\tau) \bigl(\phi \bigl(
\Phi^{-1}(\tau ) \bigr) \bigr)^{-2}h_{\FZ}^{5},
\label{(7)}
\end{equation}
where, for the selection of $h_{\FZ}$, we may rely, for instance, on
the plug-in
rule developed by \cite{ZL00}.

This rule (\ref{(7)}) can be regarded as the
combination of a plug-in strategy and a rule-of-thumb: plug-in strategy
in the selection of
$h_{\FZ}$ but rule-of-thumb for the dependence on $\tau$. It
furthermore implies that the selected $ h_{\boldsymbol{\tau}}$ has
the same $n^{-1/7}$ rate of convergence as $h_{\FZ}$ (see \cite{ZL00}).

%%%%%%%%%%%%%%%%%%%%%%%%%%%%%%%%%%%%%%%%%%
%s7 #&#
\section{A real data example}
\label{pracsec}
%%%%%%%%%%%%%%%%%%%%%%%%%%%%%%%%%%%%%%%%%%
\label{simul2}

\begin{figure}

\includegraphics{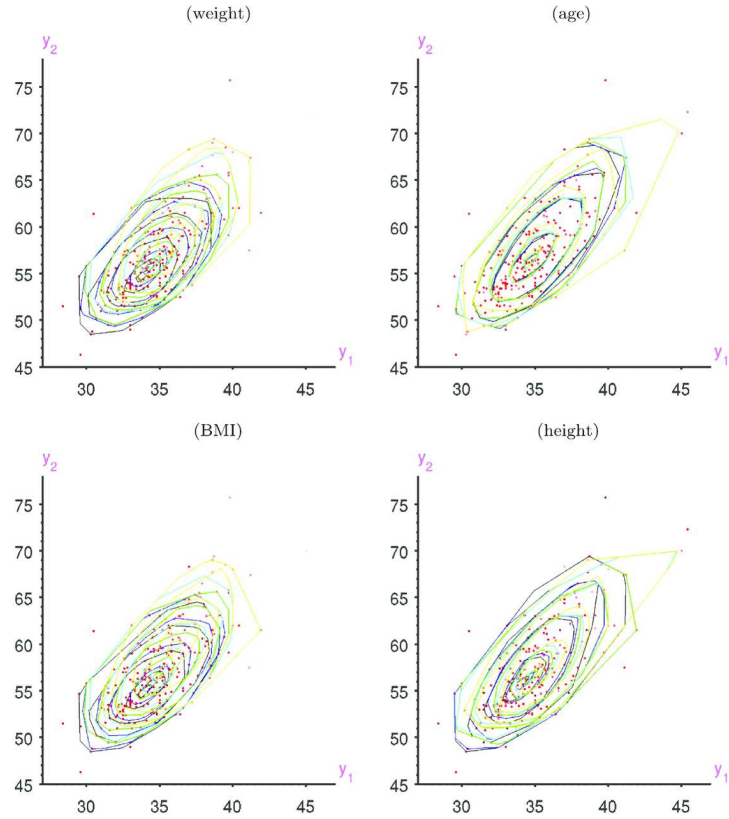}

\caption{Four empirical (local constant) regression quantile plots
from the body girth measurements dataset (women subsample; see \protect
\cite
{He03}). Throughout, the bivariate response $(Y_1,Y_2)^{\prime}$ involves
calf maximum girth ($Y_1$) and thigh maximum girth ($Y_2$), while a
single random regressor is used: weight, age, BMI, or height. The plots
are providing, for $\tau =$ 0.01, 0.03, 0.10, 0.25, and 0.40, the cuts of
the local constant regression $\tau$-quantile contours, at the
empirical $p$-quantiles of the regressors, for $p=$ 0.10 (black), 0.30
(blue), 0.50 (green), 0.70 (cyan) and 0.90 (yellow). The $n=260$ data
points are shown in red (the lighter the red color, the higher the
regressor value). The results are based on a Gaussian kernel and the
bandwidth $H = 3 \sigma_w n^{-1/5}$, where $\sigma_w$ stands for the
empirical standard deviation of the regressor (the corresponding cuts
obtained from linear regression are provided in Figure~7 of HP\v{S}).
A color version of this figure is more readable, and can be found in the on-line edition of the paper.}
\label{BGMConst}
%\vspace*{-11pt}
\end{figure}

\begin{figure}

\includegraphics{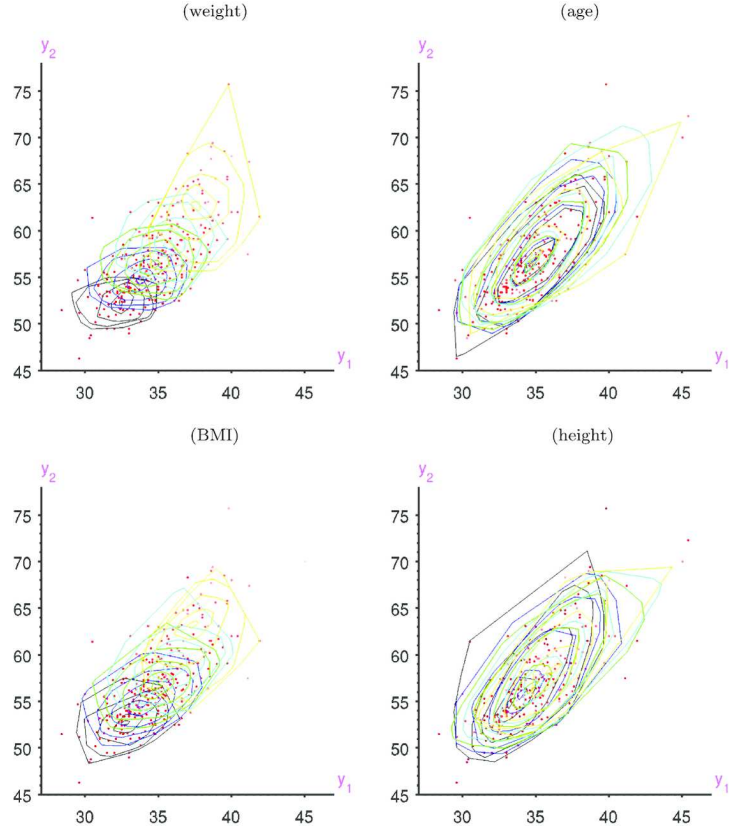}

\caption{Same quantities as in Figure \protect\ref{BGMConst}, here
obtained from the
local bilinear approach, with the same kernel and bandwidth as in
Figure \protect\ref{BGMConst} %.
%$H = 1.3\times\sigma_w n^{-1/5}$, where $\sigma_w$ stands for the
%empirical standard deviation of the regressor used).
(the computation was based on 360 equispaced directions $\mathbf{u}\in
\mathcal{S}^1$). A~color version of this figure is more readable, and can be found in the on-line edition of the paper.}
\label{BGMLin}
\vspace*{-6pt}
\end{figure}

In order to illustrate the data-analytic power of the proposed method,
we consider the ``body girth measurement'' dataset from \cite{He03},
that was already investigated in HP\v{S}. The dataset consists of
joint measurements of nine skeletal and twelve body girth dimensions,
along with weight, height, and age, in a group of 247 young men and 260
young women. As in HP\v{S}, we discard the male observations, we
restrict to the calf maximum girth ($Y_1$) and the thigh maximum girth
($Y_2$) for the response, and use a single random regressor $W$
(weight, height, age, or BMI). Figures~\ref{BGMConst} and \ref
{BGMLin} provide cuts -- for the same $w$- and $\tau$-values as in HP\v
{S} -- obtained from the proposed local constant and local bilinear
approaches, respectively.
%
%f5 #&#
%
%f6 #&#

%%%%%%%%%%%%%%%%%%%%%%%%%%%%%%%%%%%%%%%%%%%%%%

These cuts confirm most of the global analysis conducted in HP\v{S}
and moreover reveal some interesting new features. For instance,
%\begin{itemize}
%\item[(a)]
\begin{enumerate}[(b)]
\item[(a)] for the dependence on weight, the local bilinear approach confirms
the positive trend in location, the increase in dispersion, and the
evolution of ``principal directions'' (as weight increases, the first
``principal direction'' rotates from horizontal to vertical), and it
further indicates that high weights give rise to simultaneously large
extreme values in $Y_1$ and $Y_2$. The differences, for low and high
values of the covariate (weight), between the contours resulting from
the local bilinear and local constant approaches illustrate the
sensitivity of the latter to boundary effect;

\item[(b)] %\item[(b)]
for the dependence on age, the local regression quantile regions,
parallel to their global HP\v{S} counterparts, do indicate that the location
%(conditional deepest point)
and the first principal direction (along the main bisector) are
constant over age. Still as in HP\v{S}, the local approaches confirm
that the shapes of outer contours vary quite significantly with age,
indicating\vadjust{\goodbreak} an increasing (with age) simultaneous variability of both
calf and thigh girth largest values. Now, compared to HP\v{S}, the
local bilinear approach further shows that \emph{young} women present
a large simultaneous variability of both calf and thigh girth \emph
{smallest} values;

\item[(c)] %\item[(c)]
for the dependence on height, the local methods confirm the regression
effect specific to inner contours. The local bilinear approach further
shows that there is also a regression effect for outer contours that,
as height increases, get more widespread in the direction $\mathbf{u}$
(corresponding to simultaneously large values of both responses).
\end{enumerate}

Limited as it is, this short application demonstrates how the local
quantile regression analysis proposed here complements and refines the
findings obtained from the global approach introduced in HP\v{S} by
revealing the possible non-linear, heteroskedastic, skewness \ldots
features of the distributions of $\mathbf Y$ conditional on $\mathbf
{W} = \mathbf{w}$. We refer to \cite{Metal11} for a further
application, in the context of bivariate growth charts.\vadjust{\goodbreak}

We conclude this section with a brief discussion of the computational
aspects of the proposed methods. In principle, any quantile
regression/linear programming/convex optimization solver can be used
for that purpose. The exact local constant quantile/depth contours can
be computed for any $w_0$ via a weighted version of the HP\v{S}
algorithm -- see Paindaveine and {\v{S}}iman \cite{PS11b} for a detailed
description of its Matlab
implementation and its computation cost. The local bilinear contours,
for given $w_0$, are determined by considering a fixed number $M$ of
directions; their computation then is as demanding as $M$ times the
standard simple-output quantile regression with the same number of
regressors; see Koenker \cite{K05} for computational and
algorithmic details.

%%%%%%%%%%%%%%%%%%%%%%%%%%%%%%%%%%%%%%%%%%
%s8 #&#
\section{Conclusion}
\label{conclsec}
%%%%%%%%%%%%%%%%%%%%%%%%%%%%%%%%%%%%%%%%%%

In this paper, we propose a definition of regression depth as the
conditional depth of an $m$-dimensional response conditional on a
$p$-dimensional covariate. We also propose local constant and local
bilinear methods for the estimation of conditional depth contours, and
establish the consistency and asymptotic normality of the estimators.
As a descriptive tool, the resulting contours provide a powerful
data-analytic tool, while our asymptotic results guarantee that, for
$n$ large enough, those contours are able to detect any
covariate-dependent feature of the conditional distributions of the
response. An important domain of application for such methods is in the
analysis of multiple output growth charts, where current practice is
essentially restricted to a marginal approach that neglects all
information related to \textit{joint} conditional features.

%%%%%%%%%%%%%%%%%

%%%%%%%%%%%%%%%%%

\begin{appendix}

%s9 #&#
\section*{Appendix: Proofs of asymptotic results}\label{app}

We actually restrict to the local bilinear case (proofs for the local
constant case are entirely similar). The proofs rely on several lemmas,
and require some further notation.

Referring to (\ref{paramlin}) and (\ref{Zudi1.1a}), define
\[
\thetab^\ell = \vec \lleft( %
\begin{array} {c@{
\quad}c} a_{{\boldsymbol{\tau}};\mathbf{w}_0} & \mathbf{c}_{{\boldsymbol
{\tau}};\mathbf
{w}_0}^{\prime}
\\
\dot\mathbf{a}_{{\boldsymbol{\tau}};\mathbf{w}_0} & \dot\mathbf {c}_{{\boldsymbol{\tau}};\mathbf{w}_0}^{\prime}
\end{array} %
\rright) =: \vec \lleft( %
\begin{array} {c} { \boldsymbol{\varpi}}_{\mathbf{w}_0}^{\prime}
\\
\dot{\boldsymbol{\varpi}}_{\mathbf{w}_0}^{\prime} \end{array} %
\rright) \quad \mbox{and}\quad %
%and
%% \mbox{ and } $$
\hat{
\thetab}^{\ell(n)} = \vec \lleft( %
\begin{array} {c@{\quad}c}
\hat{a} & \hat{\mathbf{c}}^{\prime}
\\
\hat{\dot\mathbf{a}} & \hat{\dot\mathbf{c}}^{\prime} \end{array}
\rright) =: \vec \lleft( %
\begin{array} {c}
\widehat{ \boldsymbol{\varpi}}_{{\mathbf{w}_0}}^{\prime}
\\
\widehat{\dot{\boldsymbol{\varpi}}}_{\mathbf{w}_0} ^{\prime} \end{array}
\rright). %
\]
Denote by ${\boldsymbol{\varpi}}_1 = (a_1, \mathbf{c}_1^{\prime
})^{\prime}$ and
$\widetilde{\boldsymbol{\varpi}}_1= (\tilde{a}_1, \tilde\mathbf
{c}_1^{\prime}
)^{\prime}$ two arbitrary vectors of $\R^m$, by ${\boldsymbol{\varpi
}}_2= (\mathbf
{a}_2, \mathbf{c}_2^{\prime})^{\prime}$ and $\widetilde{\boldsymbol
{\varpi}}_2=
(\tilde\mathbf{a}_2, \tilde\mathbf{c}_2^{\prime})^{\prime}$ two arbitrary
$m\times(p-1)$ matrices.
Let $H_n:=\sqrt{nh_n^{p-1}}$ and put
%
%e9.1 #&#
\begin{eqnarray}
\label{defvarphin} {\boldsymbol{\varphi}}^{(n)}&:=& H_n
\mathbf{M}_h^\ell \vec \lleft( %
\begin{array} {c} ( \widehat{{\boldsymbol{\varpi}}}_{\mathbf{w}_0} -{{
\boldsymbol { \varpi}}}_{\mathbf
{w}_0})^{\prime}
\nonumber
\\
( \widehat{\dot{{\boldsymbol{\varpi}}}}_{\mathbf{w}_0}-{\dot {{\boldsymbol{
\varpi} }}}_{\mathbf{w}_0} )^{\prime}
\end{array}
\rright)
,
\\
{\boldsymbol{\varphi}} &:=& H_n \mathbf{M}_h^\ell
\vec \lleft( %
\begin{array} {c} ( {{\boldsymbol{
\varpi}}}_1 -{{\boldsymbol{ \varpi}}}_{\mathbf
{w}_0})^{\prime}
\\
( {{{\boldsymbol{\varpi}}}}_2-{\dot{{\boldsymbol{\varpi
}}}}_{\mathbf{w}_0} )^{\prime} \end{array} %
\rright),
\\
\tilde{\boldsymbol{\varphi}}&:=& H_n \mathbf{M}_h^\ell
\vec \lleft( %
\begin{array} {c} ( \widetilde{{
\boldsymbol{ \varpi}}}_1 -{{\boldsymbol{\varpi }}}_{\mathbf{w}_0})^{\prime}
\\
( \widetilde{{{\boldsymbol{\varpi}}}}_2-{\dot{{\boldsymbol{\varpi
}}}}_{\mathbf{w}_0} )^{\prime} \end{array} %
\rright) ,
\nonumber
\end{eqnarray}
and note that
${\boldsymbol{\varphi}}^{(n)}=\sqrt{n h^{p-1}_n}
\mathbf{M}_h^\ell
( \hat{\boldsymbol{\theta}}^{\ell(n)}-{\boldsymbol{\theta
}}^\ell ) $. Define
$\mathbf{W}_{hi}:=(\mathbf{W}_i-\mathbf{w}_0)/h_n$,
$K_{hi}:=K(\mathbf{W}_{hi}) $ and
$
\XX_{hi\mathbf{u}}^{\ell}
:=
(\mathbf{M}_h^\ell)^{-1}\XX_{i\mathbf{u}}^\ell
=
(1, \mathbf{Y}_{i\mathbf{u}}^{\perp\prime})^{\prime}
\otimes
(1,\mathbf{W}_{hi}^{\prime})^{\prime}
$.

Let ${Z}_{i\mathbf{u}}^\ell={Z}_{i\mathbf{u}}^\ell({\boldsymbol
{\theta}
}^\ell):={Y}_{i\mathbf{u}} -{\boldsymbol{\theta}}^{\ell\prime}\XX
_{i\mathbf{u}}^\ell
$ as in Theorem~\ref{bahadur}, and define
\begin{eqnarray*}
T_{ni}&:=& h_n\dot\mathbf{a}_{{\boldsymbol{\tau}};\mathbf{w}_0}^{\prime
}
\mathbf{W}_{hi} + h_n (\vec \dot{ \mathbf{c}}_{{\boldsymbol{\tau}};\mathbf
{w}_0})^{\prime}
\bigl(\mathbf{Y}_{i\mathbf{u}}^{\perp} \otimes \mathbf{W}_{hi}
\bigr), %
\\
Z^{*}_{ni}({\boldsymbol{\varphi}}) &:=&
{Z}_{i\mathbf{u}}^\ell- H_n^{-1} {\boldsymbol{
\varphi}}^{\prime
}\XX_{hi\mathbf
{u}}^{\ell}\quad \mbox{and}\quad
U_{ni}=U_{ni}({\boldsymbol{\varphi}}) := T_{ni}+
H_n^{-1} {\boldsymbol{\varphi}}^{\prime}
\XX_{hi\mathbf
{u}}^{\ell} %
\end{eqnarray*}
(note that the latter two quantities depend on the choice of
${\boldsymbol{\varpi}}_1$ and ${\boldsymbol{\varpi}}_2$). The
following identities will be
useful in the sequel:
%
%e9.2 #&#
\begin{eqnarray}
\label{B.1a}{Z}_{i\mathbf{u}}^\ell &=& Y_{i\mathbf{u}} - \bigl(
a_{{\boldsymbol{\tau}};\mathbf{w}_0} + \mathbf{c}_{{\boldsymbol{\tau}};\mathbf{w}_0}^{\prime}\mathbf
{Y}^\perp _{i\mathbf{u}} \bigr) -T_{ni},
\\
\label{B.1b}Z^{*}_{ni}({\boldsymbol{\varphi}}) &=&
Y_{i\mathbf{u}} - \bigl(a_{{\boldsymbol{\tau}};\mathbf{w}_0} + \mathbf{c}_{{\boldsymbol
{\tau}};\mathbf
{w}_0}^{\prime}
\mathbf{Y}^\perp_{i\mathbf{u}} \bigr) -U_{ni}({\boldsymbol{
\varphi}})
\nonumber
\\[-8pt]
\\[-8pt]
&=& Y_{i\mathbf{u}} - \bigl({\vec} ({{\boldsymbol{ \varpi}}}_1 ,
{{{\boldsymbol{\varpi}}}}_2)^{\prime
} \bigr) ^{\prime}
\XX_{i\mathbf{u}}^\ell %-a_0-a_1'(\mathbf{X}_{i}-x_0)-[\mathbf{c}_0+ \mathbf{c}_1(
%\mathbf{X}_{i}-x_0)]\pr(-{\boldsymbol{\Gamma}}\pr_\mathbf{u}
%\mathbf{Y}_i)
.
\nonumber
\end{eqnarray}
Let $C$ be a generic constant whose value may vary from line to
line. Since $K$ is a bounded density with a
bounded support, we have, whenever $K_{hi}>0$,
%
%e9.4 #&#
\begin{equation}
\|\mathbf{W}_{hi}\|\leq C \quad \mbox{and}\quad \bigl\|
\XX_{hi\mathbf{u}}^{\ell} \bigr\|\leq C \bigl(1+\bigl\|\mathbf{Y}^\perp
_{i\mathbf{u}}\bigr\| \bigr), \label{B.1c}
\end{equation}
and, when moreover $\|{\boldsymbol{\varphi}}\|\leq M$,
%
%e9.5 #&#
\begin{equation}
|T_{ni}|\leq Ch_n \bigl(1+\bigl\|\mathbf{Y}^\perp_{i\mathbf
{u}}
\bigr\| \bigr)\quad \mbox{and}\quad |U_{ni}|\leq C \bigl(h_n+H_n^{-1}
\bigr) \bigl(1+\bigl\|\mathbf {Y}^\perp_{i\mathbf{u}}\bigr\| \bigr) .
\label{B1d}
\end{equation}

It follows from the definition of $\hat{\thetab}^{\ell(n)}$ as the
argmin of (\ref{Zudi1.1.lc}) that
%
%e9.6 #&#
\begin{equation}
\label{B.1e} {\boldsymbol{\varphi}}^{(n)} = \mathop{\argmin}_{\boldsymbol{\varphi}\in\mathbb{R}^{mp}}
\sum_{i=1}^n K_{hi}
\rho_{\tau} \bigl( Z^{*}_{ni}(\boldsymbol{\varphi})
\bigr) %\stackrel{{\boldsymbol{\psi}}}{=}\mbox\mathrm{argmin}_{{\boldsymbol{
%\varphi}}\in
%%R^{1+p}}G_n({\boldsymbol{\varphi}})
.
\end{equation}
Recalling that $\psi_{\tau}(y) :={\tau}-I[ y < 0] $, define
%
%e9.7 #&#
\begin{equation}
\mathbf{V}_n({\boldsymbol{\varphi}}) := H_n^{-1}
\sum_{i=1}^n K_{hi}
\psi_{\tau} \bigl( Z^{*}_{ni}({\boldsymbol{\varphi
}}) \bigr)\XX _{hi\mathbf{u}}^{\ell}. \label{B.2}
\end{equation}
In order to prove Theorem~\ref{bahadur}, we need the following lemma.

%le9.1 #&#
\begin{Lem}
\label{LemA1}
Let $\mathbf{V}_n(\cdot)\dvtx \R^{mp}\to\R^{mp}$ be a sequence of
functions that satisfies the following two properties:
\begin{enumerate}[(ii)]
\item[(i)]
for all $\lambda\geq1$ and all ${\boldsymbol{\psi}}\in\R^{mp}$,
$-{\boldsymbol{\psi}}^{\prime}\mathbf{V}_n(\lambda{\boldsymbol
{\psi}})\geq-{\boldsymbol{\psi}}
^{\prime}\mathbf{V}_n({\boldsymbol{\psi}})$ a.s.;
\item[(ii)]
there exist a $p\times p$ positive definite
matrix $\mathbf{D}$ and a sequence of $mp$-dimensional random vectors
$\mathbf{A}_n$ satisfying $\|\mathbf{A}_n\|=\mathrm{O}_\mathrm{P}(1)$
such that, for all $M>0$,\break $\sup_{\|{\boldsymbol{\psi}}\|\leq M}\|
\mathbf
{V}_n({\boldsymbol{\psi}})+(\mathbf{G}_{\taub;\mathbf{w}_0}\otimes
\mathbf
{D}){\boldsymbol{\psi}}- \mathbf{A}_n\|=\mathrm{o}_\mathrm{P}(1)$,
where $ \mathbf{G}_{\taub;\mathbf{w}_0}$ is given in Assumption~\textup{\ref{as1}(v)}.
\end{enumerate}
Then, if ${\boldsymbol{\psi}}_n$ is such that
$\| \mathbf{V}_n({\boldsymbol{\psi}}_n)\|=\mathrm{o}_\mathrm{P}(1)$, it
holds that $\|
{\boldsymbol{\psi}}_n\|=\mathrm{O}_\mathrm{P}(1)$ and
%
%e9.8 #&#
\begin{equation}
\label{A.3} \boldsymbol{\psi}_n = (\mathbf{G}_{\taub;\mathbf{w}_0}
\otimes \mathbf{D})^{-1}\mathbf {A}_n+\mathrm{o}_\mathrm{P}(1).
\end{equation}
\end{Lem}
\begin{pf} The proof follows along the same lines as
in page 809 of \cite{KZ96}; details are left to the
reader.
\end{pf}

The proof of Theorem~\ref{bahadur} consists in checking that the
assumptions of
Lemma~\ref{LemA1} hold for $\mathbf{V}_n$ defined in (\ref{B.2}); we
use the following lemma.

%le9.2 #&#
\begin{Lem}\label{LemA2}
Under Assumptions \textup{\ref{as1}}--\textup{\ref{as3}}, for any $({\boldsymbol{\varphi}}, \tilde
{\boldsymbol{\varphi}})$ such that
$\max(\|{\boldsymbol{\varphi}}\|,\|\tilde{\boldsymbol{\varphi}}\|
)\leq M$, and
$n$ large enough,
%
%e9.9 #&#
\begin{eqnarray}
\label{first}\mathrm{E} \bigl[K_{hi} \bigl|\psi_{\tau} \bigl(
Z^{*}_{ni}({\boldsymbol {\varphi} }) \bigr)-
\psi_{\tau} \bigl( Z^{*}_{ni}(\tilde{\boldsymbol{
\varphi}}) \bigr) \bigr| \bigr] &\leq& C\mathrm{E} \bigl[K_{hi} I \bigl[\bigl|
Z^{*}_{ni}(\tilde{\boldsymbol{\varphi} })\bigr|<CH_n^{-1}
\|{\boldsymbol{\varphi}}-\tilde{\boldsymbol{\varphi }}\| \bigr] \bigr]
\nonumber
\\[-8pt]
\\[-8pt]
&\leq& Ch_n^{p-1}H_n^{-1}\|{
\boldsymbol{\varphi}}-\tilde{\boldsymbol {\varphi}}\|
\nonumber
\end{eqnarray}
and
%
%e9.10 #&#
\begin{eqnarray}
\label{second} \mathrm{E} \bigl[K_{hi}^2 \bigl|
\psi_{\tau} \bigl( Z^{*}_{ni}({\boldsymbol {
\varphi} }) \bigr)- \psi_{\tau} \bigl( Z^{*}_{ni}(
\tilde{\boldsymbol{ \varphi}}) \bigr)\bigr|^2 \bigr] &\leq& C\mathrm{E}
\bigl[K_{hi}^2 I \bigl[\bigl| Z^{*}_{ni}(
\tilde{\boldsymbol{\varphi} })\bigr|<CH_n^{-1}\|{\boldsymbol{
\varphi}}-\tilde{\boldsymbol{ \varphi }}\| \bigr] \bigr]
\nonumber
\\[-8pt]
\\[-8pt]
&\leq& Ch_n^{p-1}H_n^{-1}\|{
\boldsymbol{\varphi}}-\tilde{\boldsymbol {\varphi}}\|.
\nonumber
\end{eqnarray}
\end{Lem}

\begin{pf} The claim, in this lemma, is similar to
that of Lemma A.3 in \cite{Hetal09}, which essentially
follows from the same argument as in the time series case (cf. \cite
{Letal98}). Details, however, are
quite different. It follows from (\ref{B.1c}) that
\begin{eqnarray*}
%\lefteqn{
%
K_{hi} \bigl|\psi_{\tau} \bigl(
Z^{*}_{ni}({\boldsymbol{\varphi}}) \bigr)-\psi
_{\tau} \bigl( Z^{*}_{ni}(\tilde{\boldsymbol{
\varphi}}) \bigr)\bigr| &=& K_{hi} \bigl|I \bigl[ Z^{*}_{ni}({
\boldsymbol{\varphi}})<0 \bigr] -I \bigl[ Z^{*}_{ni}(
\tilde{ \boldsymbol{\varphi}})<0 \bigr] \bigr| %}
\\
& = & K_{hi} \bigl|I \bigl[ Z^{*}_{ni}(\tilde{
\boldsymbol{\varphi}})< H_n^{-1} ({\boldsymbol{\varphi} }-
\tilde{\boldsymbol{\varphi}})^{\prime}\XX_{hi\mathbf{u}}^{\ell
}
\bigr] -I \bigl[ Z^{*}_{ni}(\tilde{\boldsymbol{\varphi}})<0
\bigr] \bigr|
\\
& \leq& K_{hi} I \bigl[ \bigl| Z^{*}_{ni}(\tilde{
\boldsymbol{\varphi}})\bigr|< C H_n^{-1}\|{\boldsymbol{\varphi}}-
\tilde{\boldsymbol{\varphi}}\| \bigl(1+\bigl\|\mathbf {Y}^\perp_{i\mathbf{u}}
\bigr\| \bigr) \bigr] .
\end{eqnarray*}
Hence, from (\ref{B.1b}) and the mean value theorem, we obtain
\begin{eqnarray*}
&&\mathrm{E} \bigl[ K_{hi} \bigl|\psi_{\tau} \bigl(
Z^{*}_{ni}({\boldsymbol {\varphi} }) \bigr)-
\psi_{\tau} \bigl( Z^{*}_{ni}(\tilde{\boldsymbol{
\varphi}}) \bigr)\bigr| \bigr]
\\
&&\quad \leq \mathrm{E} \bigl[ K_{hi} I \bigl[ \bigl| Z^{*}_{ni}(
\tilde{\boldsymbol{\varphi} })\bigr|<CH_n^{-1}\|{\boldsymbol{
\varphi}}-\tilde{\boldsymbol{\varphi }}\| \bigl(1+\bigl\|\mathbf {Y}^\perp_{i\mathbf{u}}
\bigr\| \bigr) \bigr] \bigr]
\\
&&\quad = \mathrm{E} \bigl[ K_{hi} \mathrm{P} \bigl[ \bigl|
Z^{*}_{ni}( \tilde{\boldsymbol {\varphi}
})\bigr|<CH_n^{-1}\|{\boldsymbol{ \varphi}}-\tilde{\boldsymbol{
\varphi }}\| \bigl(1+\bigl\|\mathbf {Y}^\perp_{i\mathbf{u}} \bigr\| \bigr) |
\mathbf{Y}_{i\mathbf{u}}^\perp,\mathbf{W}_i \bigr] \bigr]
\\
&&\quad = \mathrm{E} \bigl[ K_{hi} F^{Y_{\mathbf{u}} \vert(\mathbf{Y}_{\mathbf
{u}}^\perp,\mathbf{W})}
\bigl(a_{{\boldsymbol{\tau}};\mathbf{w}_0} +\mathbf {c}_{{\boldsymbol{\tau}};\mathbf{w}_0}^{\prime}
\mathbf{Y}^\perp _{i\mathbf
{u}}+U_{ni}(\tilde{\boldsymbol{
\varphi}}) +CH_n^{-1}\|{\boldsymbol{\varphi}}-\tilde{
\boldsymbol{ \varphi}}\| \bigl(1+\bigl\|\mathbf {Y}^\perp_{i\mathbf{u}}\bigr\|
\bigr) %|\mathbf{Y}^\perp_{i\mathbf{u}},\mathbf{X}_i
\bigr) \bigr]
\\
&&\qquad {} - \mathrm{E} \bigl[ K_{hi} F^{Y_{\mathbf{u}} \vert(\mathbf{Y}_{\mathbf{u}}^\perp
,\mathbf{W})}
\bigl(a_{{\boldsymbol{\tau}};\mathbf{w}_0} +\mathbf {c}_{{\boldsymbol{\tau}};\mathbf{w}_0}^{\prime}
\mathbf{Y}^\perp _{i\mathbf
{u}}+U_{ni}(\tilde{\boldsymbol{
\varphi}})-CH_n^{-1} \|{\boldsymbol {\varphi}}-\tilde {
\boldsymbol{\varphi}}\| \bigl(1+\bigl\| \mathbf{Y}^\perp_{i\mathbf{u}}\bigr\|
\bigr) %|\mathbf{Y}^\perp_{i\mathbf{u}},\mathbf{X}_i
\bigr) \bigr]
\\
&&\quad \le \mathrm{E} \bigl[ K_{hi} \bigl(1+\bigl\|\mathbf{Y}^\perp_{i\mathbf{u}}
\bigr\| \bigr) f^{Y_{\mathbf{u}} \vert(\mathbf{Y}_{\mathbf{u}}^\perp,\mathbf
{W})} \bigl(a_{{\boldsymbol{\tau}};\mathbf{w}_0} +\mathbf {c}_{{\boldsymbol{\tau}
};\mathbf{w}_0}^{\prime}
\mathbf{Y}^\perp_{i\mathbf{u}}+ U_{ni} (\tilde{\boldsymbol{
\varphi}})
\\
&&\hphantom{\quad \le \mathrm{E} \bigl[ K_{hi} \bigl(1+\|
\mathbf{Y}^\perp_{i\mathbf{u}} \| \bigr) f^{Y_{\mathbf{u}} \vert(\mathbf{Y}_{\mathbf{u}}^\perp,\mathbf
{W})} \bigl(} {}+
\lambda CH_n^{-1}\|{\boldsymbol {\varphi}} - \tilde{
\boldsymbol{\varphi}}\| \bigl(1+\bigl\|\mathbf{Y}^\perp_{i\mathbf
{u}} \bigr\|
\bigr) \bigr) \bigr]
\\
&&\qquad {} \times2CH_n^{-1}\|{\boldsymbol{\varphi}}-
\tilde{ \boldsymbol {\varphi}}\|
\end{eqnarray*}
for some $\lambda\in(-1,1)$.
Assumptions \ref{as1}--\ref{as3},
together with (\ref{B1d}), therefore yield that, for ${\boldsymbol
{\varphi}},
\tilde{\boldsymbol{\varphi}}\in\{{\boldsymbol{\varphi}}\dvtx
\|{\boldsymbol{\varphi}}\|\leq M\}$ and $n$ large enough,\vspace*{2pt}
\begin{eqnarray*}
&&\mathrm{E} \bigl[ K_{hi} \bigl|\psi_{\tau} \bigl(
Z^{*}_{ni}({\boldsymbol{\varphi}}) \bigr)-
\psi_{\tau} \bigl( Z^{*}_{ni}(\tilde{\boldsymbol{
\varphi}}) \bigr)\bigr| \bigr]
\\
&&\quad \leq CH_n^{-1}\|{\boldsymbol{\varphi}}-\tilde{
\boldsymbol{\varphi}}\| %\
%\
%& & \times
\mathrm{E} \biggl[{
K_{hi} \int_{\mathbb{R}^{m-1}}} \bigl(1+\|\mathbf{t}
\|\bigr)f^{Y_{\mathbf{u}} \vert(\mathbf{Y}_{\mathbf
{u}}^\perp=\mathbf{t} ,\mathbf{W})} \bigl(a_{{\boldsymbol{\tau}};\mathbf{w}_0} +\mathbf{c}_{{\boldsymbol{\tau}};\mathbf{w}_0}^{\prime}
\mathbf{t} \bigr) f^{ \mathbf{Y}_{\mathbf{u}}^\perp\vert\mathbf{W}}(\mathbf{t} ) \,\mathrm{d}\mathbf{t} \biggr]
\\
&&\quad = C h_n^{p-1}H_n^{-1} \|{
\boldsymbol{\varphi}}-\tilde{\boldsymbol {\varphi}}\| f^\mathbf{W}(
\mathbf{w}_0)
\\
&&\qquad {} \times {\int_{ \mathbb{R}^{m-1} } } \bigl(1+\|\mathbf{t}\|\bigr)
f^{ Y_{\mathbf{u}} \vert( \mathbf{Y}_{\mathbf{u}}^\perp=\mathbf
{t} ,\mathbf{W}=\mathbf{w}_0)} \bigl(a_{{\boldsymbol{\tau}};\mathbf{w}_0} +\mathbf{c}_{{\boldsymbol{\tau}};\mathbf{w}_0}^{\prime}
\mathbf{t} \bigr) f^{ \mathbf{Y}_{\mathbf{u}}^\perp\vert\mathbf{W}=\mathbf{w}_0} (\mathbf{t} ) \,\mathrm{d}\mathbf{t} ,
\end{eqnarray*}
which establishes (\ref{first}); (\ref{second}) follows along similar
lines.
\end{pf}
%
%le9.3 #&#
\begin{Lem}\label{LemA3}
Under Assumptions \textup{\ref{as1}}--\textup{\ref{as3}},
we have that, as $\ny$,\vspace*{2pt}
%
%e9.11 #&#
\begin{equation}\label{B.3}
\sup_{\|{\boldsymbol{\varphi}}\|\leq M} \bigl\| \mathbf{V}_n({\boldsymbol{
\varphi}})-\mathbf{V}_n(\mathbf{0}) -\mathrm{E} \bigl[
\mathbf{V}_n({\boldsymbol{\varphi}})-\mathbf {V}_n(
\mathbf{0}) \bigr] \bigr\| =\mathrm{o}_\mathrm{P}(1).
\end{equation}
\end{Lem}
\begin{pf} The proof of this lemma is quite similar, in view of
Lemma~\ref{LemA2}, to that of Lemma A.4 in \cite{Hetal09}. Details
are therefore omitted.
\end{pf}
%
%le9.4 #&#
\begin{Lem}\label{LemA4}
Under Assumptions \textup{\ref{as1}}--\textup{\ref{as3}},
we have that, as $\ny$,\vspace*{2pt}
%
%e9.12 #&#
\begin{equation}
\sup_{\|{\boldsymbol{\varphi}}\|\leq M}\bigl\|\mathrm{E} \bigl[\mathbf {V}_n({
\boldsymbol{\varphi}})-\mathbf{V}_n(\mathbf{0}) \bigr]+(\mathbf
{G}_{\taub;\mathbf
{w}_0}\otimes\mathbf{D}){\boldsymbol{\varphi}}\bigr\|=\mathrm{o}(1),
\label{B.12}
\end{equation}
where $\mathbf{D}=f^\mathbf{W}(\mathbf{w}_0)
\diag(1, {\boldsymbol{\mu}}_2^K)$.
\end{Lem}
\begin{pf} Note that
$%\begin{align}
\mathbf{V}_{n}({\boldsymbol{\varphi}})-\mathbf{V}_{n}(\mathbf{0})
%&
=
H_{n}^{-1}
\sum_{i=1}^n
K_{hi}
[\psi_{\tau}( Z^{*}_{ni} ({\boldsymbol{\varphi}}))
-
\psi_{\tau} ( {Z}_{i\mathbf{u}}^\ell)]
\XX_{hi\mathbf{u}}^{\ell}$.
%=:
%H_{n}^{-1}
%\sum_{i=1}^n
%\mathbf{V}_{n i}({\pmb{\theta}}),
%
%$ \label{(B.5)} \end{align}
%where
%$\mathbf{V}_{n i}({\pmb{\theta}}):=(\mathbf{V}_{n i}^0({\pmb{
%\theta}}), ({ V}_{n
%i}^1({\pmb{\theta}}))\pr, (\mathbf{V}_{n i}^2({\pmb{\theta}}))\pr, ({
%V}_{n
%i}^3({\pmb{\theta}}))\pr)\pr$, with
%$${ V}_{n i}^0({\pmb{\theta}}):=[\psi_{\tau}( Z^{*}_{ni}({\pmb{
%\theta}}))-\psi_{\tau} (
%\mathcal{Y}_{1,i} ^{*})]K_{i} \mbox{and} \mathbf{V}_{n
%i}^1({\pmb{\theta}})=[\psi_{\tau}( Z^{*}_{ni}({\pmb{\theta}}))-\psi_{
%\tau} (
%\mathcal{Y}_{1,i} ^{*})] (-\mathcal{Y}_{-1, i})K_{i} ,$$
%$${ V}_{n i}^2({\pmb{\theta}}):=[\psi_{\tau}( Z^{*}_{ni}({\pmb{
%\theta}}))-\psi_{\tau} (
%\mathcal{Y}_{1,i} ^{*})]\mathbf{X}_{h i}K_{i} \mbox{and} \mathbf{V}_{n
%i}^3({\pmb{\theta}})=[\psi_{\tau}( Z^{*}_{ni}({\pmb{\theta}}))-\psi_{
%\tau} (
%\mathcal{Y}_{1,i} ^{*})](-\mathcal{Y}_{-1, i}) \mathbf{X}_{h i}K_{i}
%.$$
It follows from (\ref{B.1a}) and (\ref{B.1b}) that\vspace*{2pt}
\begin{eqnarray*}
%\lefteqn{
%
%%
\mathrm{E} \bigl[\mathbf{V}_n({
\boldsymbol{\varphi}})-\mathbf {V}_n(\mathbf{0}) \bigr] & = & n
H_n^{-1} \mathrm{E} \bigl[ K_{hi} \bigl(I \bigl[
{Z}_{i\mathbf{u}}^\ell<0 \bigr] - I \bigl[ Z^{*}_{ni}({
\boldsymbol{\varphi}})<0 \bigr] \bigr) \XX_{hi\mathbf{u}}^{\ell} \bigr]
%}
\nonumber
\\
%& &
%
& = & H_n h_n^{-(p-1)}
\mathrm{E} \bigl[ K_{hi} \bigl( F^{Y_{\mathbf{u}}\vert(\mathbf{Y}^\perp_{\mathbf{u}},\mathbf
{W})} \bigl(a_{{\boldsymbol{\tau}};\mathbf{w}_0} +
\mathbf {c}_{{\boldsymbol{\tau}};\mathbf
{w}_0}^{\prime}\mathbf{Y}^\perp_{i\mathbf{u}}+T_{ni}
\bigr)
\nonumber
\\
& &\hphantom{H_n h_n^{-(p-1)} \mathrm{E} \bigl[ K_{hi} \bigl(}
{} - F^{Y_{\mathbf{u}}\vert(\mathbf{Y}^\perp_{\mathbf{u}},\mathbf
{W})} \bigl(a_{{\boldsymbol{\tau}};\mathbf{w}_0} +\mathbf {c}_{{\boldsymbol{\tau}};\mathbf
{w}_0}^{\prime}
\mathbf{Y}^\perp_{i\mathbf{u}}+U_{ni} \bigr) \bigr)
\XX_{hi\mathbf{u}}^{\ell} \bigr] .
\end{eqnarray*}
Then, similar to the proof of Lemma~\ref{LemA2}, by the mean value theorem,
since $U_{ni}-T_{ni}=H_n^{-1}\XX_{hi\mathbf{u}}^{\ell\prime}
{\boldsymbol{\varphi}}$, there exists $\xi\in(0,1)$ such that\vspace*{2pt}
\begin{eqnarray*}
&&\sup_{\|{\boldsymbol{\varphi}}\|\leq M} \bigl\|\mathrm{E} \bigl[\mathbf{V}_n({
\boldsymbol{\varphi}})-\mathbf {V}_n(\mathbf {0}) \bigr]+(
\mathbf{G}_{\taub;\mathbf{w}_0}\otimes\mathbf {D}){\boldsymbol{\varphi}}\bigr\|
\nonumber
\\
&&\quad = \sup_{\|{\boldsymbol{\varphi}}\|\leq M} \bigl\| (\mathbf{G}_{\taub;\mathbf{w}_0}\otimes
\mathbf{D}) {\boldsymbol{\varphi}}
\nonumber
\\[1pt]
& &\qquad {} - h_n^{-(p-1)} \mathrm{E} \bigl[
K_{hi} f^{Y_{\mathbf{u}}\vert(\mathbf{Y}^\perp
_{\mathbf{u}},\mathbf{W})} \bigl(a_{{\boldsymbol{\tau}};\mathbf{w}_0} +\mathbf
{c}_{{\boldsymbol{\tau}};\mathbf{w}_0}^{\prime} \mathbf{Y}^\perp _{i\mathbf
{u}}+T_{ni}
+\xi H_n^{-1} \XX_{hi\mathbf{u}}^{\ell\prime} {
\boldsymbol{ \varphi}} \bigr) \XX_{hi\mathbf{u}}^{\ell}
\XX_{hi\mathbf{u}}^{\ell\prime} {\boldsymbol{\varphi}} \bigr] %\mathbf{X}^{\mathbf{u}\prime}_{hi}
\bigr\|
\nonumber
\\[1.5pt]
&&\quad = \sup_{\|{\boldsymbol{\varphi}}\|\leq M}\bigl \| \bigl\{ (\mathbf{G}_{\taub;\mathbf{w}_0}
\otimes\mathbf{D}) - h_n^{-(p-1)} \mathrm{E} \bigl[
K_{hi} f^{Y_{\mathbf{u}}\vert(\mathbf{Y}_{\mathbf{u}}^\perp,\mathbf
{W})} \bigl(a_{{\boldsymbol{\tau}};\mathbf{w}_0} +\mathbf
{c}_{{\boldsymbol{\tau}};\mathbf
{w}_0}^{\prime}\mathbf{Y}^\perp_{i\mathbf{u}}
\bigr) \XX_{hi\mathbf{u}}^{\ell}\XX_{hi\mathbf{u}}^{\ell\prime} \bigr]
\bigr\} {\boldsymbol{\varphi}}
\nonumber
\\[1.5pt]
& &\hphantom{\quad = \sup_{\|{\boldsymbol{\varphi}}\|\leq M} \|} {} -h_n^{-(p-1)}
\mathrm{E} \bigl[ K_{hi} \bigl( f^{Y_{\mathbf{u}}\vert(\mathbf{Y}^\perp_{\mathbf{u}},\mathbf
{W})} \bigl(a_{{\boldsymbol{\tau}};\mathbf{w}_0} +
\mathbf {c}_{{\boldsymbol{\tau}};\mathbf
{w}_0}^{\prime} \mathbf{Y}^\perp_{i\mathbf{u}}+T_{ni}
+\xi H_n^{-1} \XX _{hi\mathbf{u}}^{\ell\prime} {
\boldsymbol{ \varphi}}\bigr)
\\[1.5pt]
& &\hphantom{\hphantom{\quad = \sup_{\|{\boldsymbol{\varphi}}\|\leq M} \|} {}
-h_n^{-(p-1)} \mathrm{E} \bigl[ K_{hi} \bigl( } {} - f^{Y_{\mathbf{u}}\vert(\mathbf{Y}_{\mathbf{u}}^\perp,\mathbf
{W})}
\bigl(a_{{\boldsymbol{\tau}};\mathbf{w}_0} +\mathbf {c}_{{\boldsymbol{\tau}};\mathbf
{w}_0}^{\prime}
\mathbf{Y}^\perp_{i\mathbf{u}} \bigr) \bigr) \XX_{hi\mathbf{u}}^{\ell}
\XX_{hi\mathbf{u}}^{\ell\prime} {\boldsymbol{\varphi}} \bigr] %\mathbf{X}^{\mathbf{u}\prime}_{hi}
\bigr\|
\nonumber
%&=&
% \sup_{\|{\boldsymbol{\varphi}}\|\leq M}
% \|-h_n^{-(p-1)}
% \mathrm{E}\big[ K_{hi} f^{Y_{\mathbf{u}}\vert(\mathbf{Y}^\perp_{
%\mathbf{u}},\mathbf{W})}(a_{{\boldsymbol{\tau}};\mathbf{w}_0} +
%\mathbf{c}_{{\boldsymbol{\tau}};\mathbf{w}_0}\pr\mathbf{Y}^\perp_{i
%\mathbf{u}}+T_{ni} + \xi
%H_n^{-1} {\boldsymbol{\varphi}}\pr\XX_{hi\mathbf{u}}^{\ell} )
% \nonumber
% \\
%& &
% -f^{Y_{\mathbf{u}}\vert(\mathbf{Y}_{\mathbf{u}}^\perp,
%\mathbf{W})}(a_{{\boldsymbol{\tau}};\mathbf{w}_0} +\mathbf{c}_{{
%\boldsymbol{\tau}};
%\mathbf{w}_0}\pr\mathbf{Y}^\perp_{i\mathbf{u}}) \big] {\boldsymbol{
%\varphi}}
%\pr\XX_{hi\mathbf{u}}^{\ell}
% %\mathbf{X}^{\mathbf{u}\prime}_{hi}
% \nonumber
%\\
\\[1.5pt]
&&\quad \leq C\bigl \| (\mathbf{G}_{\taub;\mathbf{w}_0}\otimes\mathbf{D}) -
h_n^{-(p-1)} \mathrm{E} \bigl[ K_{hi}
f^{Y_{\mathbf{u}}\vert(\mathbf{Y}^\perp_{\mathbf{u}},\mathbf
{W})} \bigl(a_{{\boldsymbol{\tau}};\mathbf{w}_0} +\mathbf {c}_{{\boldsymbol{\tau}};\mathbf
{w}_0}^{\prime}
\mathbf{Y}^\perp_{i\mathbf{u}} \bigr) \XX_{hi\mathbf{u}}^{\ell}
\XX_{hi\mathbf{u}}^{\ell\prime} %\mathbf{X}^{\mathbf{u}\prime}_{hi}
\bigr] \bigr\|
\\[1.5pt]
& &\qquad {} + C \sup_{\|{\boldsymbol{\varphi}}\|\leq M} h_n^{-(p-1)}
\mathrm{E} \bigl[ K_{hi} \bigl|f^{Y_{\mathbf{u}}\vert(\mathbf{Y}^\perp_{\mathbf{u}},\mathbf
{W})} \bigl(a_{{\boldsymbol{\tau}};\mathbf{w}_0} +\mathbf
{c}_{{\boldsymbol{\tau}};\mathbf
{w}_0}^{\prime}\mathbf{Y}^\perp_{i\mathbf{u}}+T_{ni}
+\xi H_n^{-1} \XX_{hi\mathbf{u}}^{\ell\prime} {
\boldsymbol{\varphi}} \bigr)
\nonumber
\\[1.5pt]
& &\hphantom{\qquad {} + C \sup_{\|{\boldsymbol{\varphi}}\|\leq M} h_n^{-(p-1)}
\mathrm{E} \bigl[ K_{hi} \bigl|} {} -f^{Y_{\mathbf{u}}\vert(\mathbf{Y}_{\mathbf{u}}^\perp,\mathbf
{W})} \bigl(a_{{\boldsymbol{\tau}};\mathbf{w}_0} +\mathbf
{c}_{{\boldsymbol{\tau}};\mathbf
{w}_0}^{\prime} \mathbf{Y}^\perp_{i\mathbf{u}}
\bigr)\bigr| \bigl\|\XX_{hi\mathbf{u}}^{\ell} \XX_{hi\mathbf{u}}^{\ell\prime}\bigr\|
\bigr] %\mathbf{X}^{\mathbf{u}\prime}_{hi} \|
% \\
%&= &
= \mathrm{o}(1) ,
\nonumber
\end{eqnarray*}
where we used Assumptions \ref{as1} and \ref{as2}, together with (\ref{B1d}).
\end{pf}
%
%le9.5 #&#
\begin{Lem}\label{LemA5}
Let Assumptions \textup{\ref{as2}} and \textup{\ref{as3}} hold. Then the random vector
${\boldsymbol{\varphi}}^{(n)}$
%be the minimizer of $\mathbf{V}_n({\boldsymbol{\varphi}})$
defined in (\textup{\ref{defvarphin}}) satisfies
$
\|\mathbf{V}_n({\boldsymbol{\varphi}}^{(n)})\|
=\mathrm{o}_\mathrm{P}(1)
%\leq
%H_n^{-1}
%\text{\mathrm{dim}}(\XX_{hi\mathbf{u}}^{\ell})\max_
%{1\leq i\leq n}\|K_{hi} \XX_{hi\mathbf{u}}^\ell\|
$.
\end{Lem}
\begin{pf} The proof follows from a argument similar to that of
Lemma A.2 on page 836 of \cite{RC80}.
\end{pf}
%
%le9.6 #&#
\begin{Lem}
\label{LemA6}
Under Assumptions \textup{\ref{as1}}--\textup{\ref{as3}}, for any $\mathbf{d}\in\mathbb{R}^{mp}$,\vspace*{2pt}
\begin{eqnarray*}
&&\lim_{\ny}\mathrm{E} \bigl[ \bigl\{ \mathbf{d}^{\prime}
\bigl(\mathbf{V}_n(\mathbf{0})-\mathrm{E} \bigl[\mathbf{V}_n(
\mathbf{0}) \bigr] \bigr) \bigr\}^2 \bigr]
\nonumber
\\[1.5pt]
& &\quad = \tau(1-\tau) f^\mathbf{W}(\mathbf{w}_0) \int
_{\mathbb{R}^{p-1}} \int_{\mathbb{R}^{m-1}} \bigl( \bigl[ \bigl(1,
\mathbf{t}' \bigr)\otimes \bigl(1,\mathbf{w}' \bigr)
\bigr] \mathbf{d} \bigr)^2 f^{\mathbf{Y}_\mathbf{u}^{\perp}\vert\mathbf{W}=\mathbf
{w}_0}(\mathbf{t})
K^2(\mathbf{w}) \,\mathrm{d}\mathbf{t}\, \mathrm{d}\mathbf{w}.
\nonumber
\end{eqnarray*}
\end{Lem}
\begin{pf} Set
$\tilde{v}_i
=
K_{hi}
\psi_{\tau}( {Z}_{i\mathbf{u}}^\ell)
\mathbf{d}^{\prime}\XX_{hi\mathbf{u}}^{\ell}
=
K_{hi}
\psi_{\tau}( {Z}_{i\mathbf{u}}^\ell)
[(1,\mathbf{Y}_{i\mathbf{u}}^{\perp\prime})\otimes(1,\mathbf{W}_{hi}')]
\mathbf{d}
$.
A simple calculation yields\vspace*{2pt}
%
%e9.13 #&#
\begin{equation}
\label{B.14} \mathrm{E} \bigl[ \bigl\{ \mathbf{d}^{\prime} \bigl(
\mathbf{V}_n( \mathbf{0})-\mathrm{E} \bigl[\mathbf{V}_n(
\mathbf{0}) \bigr] \bigr) \bigr\}^2 \bigr] = H_n^{-2}n
\Var[ \tilde{v}_1] = h_n^{-(p-1)} \Var[
\tilde{v}_1].
\end{equation}
%

%%%%%%%%%%%%%%%%%%%%%%
Note that, for $k=1,2$,\vspace*{2pt}
\begin{eqnarray*}
&&\lim_{\ny} h_n^{-(p-1)}
\mathrm{E} \bigl[ K_{h1}^k I \bigl[ {Z}_{1\mathbf{u}}^\ell<0
\bigr] \bigl( \mathbf{d}^{\prime}\XX_{h1\mathbf{u}}^{\ell}
\bigr)^k \bigr]
\nonumber
\\[1pt]
&&\quad = \lim_{\ny} h_n^{-(p-1)}
\mathrm{E} \bigl[ K_{h1}^k F^{\mathbf{Y}_{\mathbf{u}}\vert(\mathbf{Y}_{\mathbf{u}}^\perp
,\mathbf{W})}
\bigl(a_{{\boldsymbol{\tau}};\mathbf{w}_0} + \mathbf{c}_{{\boldsymbol
{\tau}};\mathbf
{w}_0}^{\prime}
\mathbf{Y}_{1\mathbf{u}}^\perp+T_{n1} \bigr) \bigl(
\mathbf{d}^{\prime}\XX_{h1\mathbf{u}}^{\ell} \bigr)^k
\bigr]
\nonumber
\\[1pt]
&&\quad = \tau f^\mathbf{W}(\mathbf{w}_0) \int
_{\mathbb{R}^{p-1}} \int_{\mathbb{R}^{m-1}} K^k(
\mathbf{w}) \bigl( \bigl[ \bigl(1,\mathbf{t}' \bigr)\otimes \bigl(1,
\mathbf{w}' \bigr) \bigr] \mathbf{d} \bigr)^k
%\nonumber\\
%&
%\times
f^{\mathbf{Y}^{\perp}_\mathbf{u}\vert\mathbf{W} = \mathbf
{w}_0}(\mathbf{t} )\, \mathrm{d}
\mathbf{t} \,\mathrm{d}\mathbf{w} ,
\nonumber
\end{eqnarray*}
%
%%%%%%%%%%%%%%%%%%%%%
%and
% \begin{eqnarray*}
%\lefteqn{
%%
%\lim_{\ny}
%\mathrm{E}
%\Big[
%K_{h1}
%I[ {Z}_{1\mathbf{u}}^\ell<0]
%(
%\mathbf{d}\pr\XX_{h1\mathbf{u}}^{\ell}
%)
%\Big]
%}
%\nonumber
%\\
%& =&
% \tau f^\mathbf{W}(\mathbf{w}_0)
%\int_{\mathbb{R}^{p-1}}
%\int_{\mathbb{R}^{m-1}}
%(
%[(1,\mathbf{t}')\otimes(1,\mathbf{w}')] \mathbf{d}
%)
% f^{\mathbf{Y}^{\perp}_\mathbf{u}\vert\mathbf{W}=\mathbf{w}_0}(
%\mathbf{t}) K(\mathbf{w}) \mathrm{d}\mathbf{t}
% \mathrm{d}\mathbf{w},
%\nonumber
%\end{eqnarray*}
%%%%%%%%%%%%%%%%%%%%%
which leads to
\begin{eqnarray*}
&&\lim_{\ny} h_n^{-(p-1)} \mathrm{E} [
\tilde{v}_1 ]
\\
&&\quad = \lim_{\ny} h_n^{-(p-1)}
\mathrm{E} \bigl[ K_{h1} \bigl( \tau - I \bigl[ {Z}_{1\mathbf{u}}^\ell<0
\bigr] \bigr) \bigl( \mathbf{d}^{\prime}\XX_{h1\mathbf{u}}^{\ell}
\bigr) \bigr]
\\
&&\quad = (\tau-\tau) f^\mathbf{W}(\mathbf{w}_0) \int
_{\mathbb{R}^{p-1}} \int_{\mathbb{R}^{m-1}} K(\mathbf{w}) \bigl(
\bigl[ \bigl(1,\mathbf{t}' \bigr)\otimes \bigl(1,
\mathbf{w}' \bigr) \bigr] \mathbf{d} \bigr) f^{\mathbf{Y}^{\perp}_\mathbf{u}\vert\mathbf{W}=\mathbf
{w}_0}(
\mathbf{t}) \,\mathrm{d}\mathbf{t} \,\mathrm{d}\mathbf{w} =0
\nonumber
\end{eqnarray*}
%
%%%%%%%%%%%%%%%%%%%%%%
and
\begin{eqnarray*}
&&\lim_{\ny} h_n^{-(p-1)} \mathrm{E} \bigl[
\tilde{v}_1^2 \bigr]
\\
&&\quad = \lim_{\ny} h_n^{-(p-1)}
\mathrm{E} \bigl[ K_{h1}^2 \bigl( \tau^2 - 2
\tau I \bigl[Z_{1\mathbf{u}}^\ell<0 \bigr] + I \bigl[Z_{1\mathbf{u}}^\ell<0
\bigr] \bigr) \bigl( \mathbf{d}^{\prime}\XX_{h1\mathbf{u}}^{\ell}
\bigr)^2 \bigr]
\nonumber
\\
&&\quad = \tau(1-\tau) f^\mathbf{W}(\mathbf{w}_0) \int
_{\mathbb{R}^{p-1}} \int_{\mathbb{R}^{m-1}} K^2(
\mathbf{w}) \bigl( \bigl[ \bigl(1,\mathbf{t}' \bigr)\otimes \bigl(1,
\mathbf{w}' \bigr) \bigr] \mathbf{d} \bigr)^2
f^{\mathbf{Y}^{\perp}_\mathbf{u}|\mathbf{W}=\mathbf{w}_0}(\mathbf {t}) \,\mathrm{d}\mathbf{t} \,\mathrm{d}\mathbf{w} .
\nonumber
\end{eqnarray*}
%
%%%%%%%%%%%%%%%%%%%%%%
Therefore,
\begin{eqnarray*}
&&\lim_{\ny} h_n^{-(p-1)} \Var[
\tilde{v}_1]
\\
&&\quad = \lim_{\ny} \bigl( h_n^{-(p-1)}
\mathrm{E} \bigl[\tilde{v}_1^2 \bigr] %\right)
-
%\left(
h_n^{-(p-1)} \bigl(\mathrm{E}[\tilde{v}_1]
\bigr)^2 \bigr)
\nonumber
\\
&&\quad = \tau(1-\tau) f^\mathbf{W}(\mathbf{w}_0) \int
_{\mathbb{R}^{p-1}} \int_{\mathbb{R}^{m-1}} K^2(
\mathbf{w}) \bigl( \bigl[ \bigl(1,\mathbf{t}' \bigr)\otimes \bigl(1,
\mathbf{w}' \bigr) \bigr] \mathbf{d} \bigr)^2
f^{\mathbf{Y}^{\perp}_\mathbf{u}\vert\mathbf{W}=\mathbf
{w}_0}(\mathbf{t})\, \mathrm{d}\mathbf{t}\, \mathrm{d}\mathbf{w},
\end{eqnarray*}
which, together with (\ref{B.14}), establishes the result.
\end{pf}
%
%%%
%
\begin{pf*}{Proof of Theorem~\ref{bahadur}} The proof consists in
checking that the conditions
of Lemma~\ref{LemA1} are satisfied. Lemmas \ref{LemA3} and \ref
{LemA4} entail that Lemma~\ref{LemA1}(ii) holds,
with
$
\mathbf{D}=f^{\mathbf{W}}(\mathbf{w}_0)
\diag(1, {\boldsymbol{\mu}}_2^K)
$
(yielding
$(\mathbf{G}_{\taub;\mathbf{w}_0}\otimes\mathbf
{D})^{-1}={\boldsymbol{\eta}}^\ell_{\taub;\mathbf{w}_0}$) and
$
\mathbf{A}_n=
\mathbf{V}_n(\mathbf{0})=
H_n^{-1}
\sum_{i=1}^n
K_{hi}
\psi_{\tau}( {Z}_{i\mathbf{u}}^\ell)
%(\mathbf{M}_h^\ell)^{-1}
\XX^\ell_{hi\mathbf{u}}$,
which, by Lemma~\ref{LemA6}, is $\mathrm{O}_\mathrm{P}(1)$.
%As for condition (i), take $\mathbf{A}_n=\mathbf{V}_n(\mathbf{0})$.
%Since
As for Lemma~\ref{LemA1}(ii), the fact that
\[
\lambda \mapsto -{\boldsymbol{\varphi}}^{\prime}
\mathbf{V}_n(\lambda{\boldsymbol {\varphi}}) = H_n^{-1}
\sum_{i=1}^n K_{hi}
\psi_{\tau} \bigl(Z_{i\mathbf{u}}^\ell-\lambda
H_n^{-1} {\boldsymbol {\varphi} }^{\prime}
\XX_{hi\mathbf{u}}^{\ell} \bigr) \bigl(-{\boldsymbol{
\varphi}}^{\prime} \XX_{hi\mathbf{u}}^{\ell} \bigr) %
\]
is non-decreasing directly follows from the fact $y\mapsto\psi_{\tau
}(y)$ is non-decreasing. Since (Lemma~\ref{LemA5} and Assumptions
\ref{as2} and \ref{as3}) $\|\mathbf{V}_n({\boldsymbol{\varphi}}^{(n)})\|$ is
$\mathrm{o}_\mathrm
{P}(1)$, Lemma~\ref{LemA1} applies, which concludes the proof.
\end{pf*}
%
%%%
%
\begin{pf*}{Proof of Theorem~\ref{a.normality}}
On the basis of the Bahadur representation of
Theorem~\ref{bahadur}, the asymptotic normality of $\hat{\thetab}
^{\ell(n)}$ follows exactly as in the
corresponding proofs for usual nonparametric regression in the
i.i.d. case (see, e.g., \cite{FG96}), yielding the
asymptotic normality with the bias (i.e., the expectation)
of the first term on the right-hand side of (\ref{bahadur.1}) as
\begin{eqnarray*}
&&\mathrm{E} \Biggl[ \frac{\eta^\ell_{{\boldsymbol{\tau}};\mathbf{w}_0}}{\sqrt{n h_n^{p-1}}} \sum
_{i=1}^n K_{h1} \psi_{\tau}
\bigl({Z}_{i\mathbf{u}}^\ell \bigr) \XX_{hi\mathbf{u}}^{\ell}
\Biggr]
\\
&&\quad = %\\[2mm]
%&=&
\frac{ \eta^\ell_{{\boldsymbol{\tau}};\mathbf{w}_0}}{\sqrt{n h_n^{p-1}}} n \mathrm{E} \bigl[
K_{h1} \psi_{\tau} \bigl({Z}_{1\mathbf{u}}^\ell
\bigr) \XX_{h1\mathbf{u}}^{\ell} \bigr]
\\
&&\quad = \eta^\ell_{{\boldsymbol{\tau}};\mathbf{w}_0}{\sqrt{n
h_n^{p-1}}} h_n^{-(p-1)} \mathrm{E} \bigl[
K_{h1} \bigl( F^{ Y_{\mathbf{u}}\vert( \mathbf{Y}_{\mathbf{u}}^\perp,\mathbf
{W})} \bigl(a_{{\boldsymbol{\tau}};\mathbf{W}} +
\mathbf{c}_{{\boldsymbol
{\tau}};\mathbf
{W}}^{\prime} \mathbf{Y}^\perp_{1\mathbf{u}}
\bigr)
\\
&&\hphantom{\quad = \eta^\ell_{{\boldsymbol{\tau}};\mathbf{w}_0}{\sqrt{n
h_n^{p-1}}} h_n^{-(p-1)} \mathrm{E}
\bigl[
K_{h1} \bigl(} {}- F^{ Y_{\mathbf{u}}|( \mathbf{Y}^\perp_{\mathbf{u}},\mathbf
{W})} \bigl(a_{{\boldsymbol{\tau}};\mathbf{w}_0} +\mathbf {c}_{{\boldsymbol{\tau}};\mathbf
{w}_0}^{\prime}
\mathbf{Y}^\perp_{1\mathbf{u}}+T_{n1} \bigr) \bigr)
\XX_{h1\mathbf{u}}^{\ell} \bigr]
\\
&&\quad = \sqrt{n h_n^{p-1}} \biggl(
\frac{h_n ^2}{2} \mathbf{B}^\ell_{\mathbf{w}_0} +\mathrm{o}
\bigl(h_n^2 \bigr) \biggr),
\end{eqnarray*}
where the last equality is derived from a first-order
Taylor expansion of $y\mapsto F^{ Y_\mathbf{u}\vert( \mathbf
{Y}_\mathbf{u}^{\perp},\mathbf{X})}(y)$ and a second-order Taylor
expansion of $\mathbf{w}\mapsto(a_{{\boldsymbol{\tau}};\mathbf
{w}},\mathbf
{c}_{{\boldsymbol{\tau}};\mathbf{w}}^{\prime})^{\prime}$ at
$\mathbf{w}=\mathbf{w}_0$
(these expansions exist in view of
Assumptions \ref{as1} and \ref{as4}). The $\mathrm{o}(h_n^2)$ term is taken care of by
Assumption~\ref{as3}. The asymptotic variance of the theorem readily follows
from Lemma~\ref{LemA6}. Details
are omitted.
\end{pf*}

\end{appendix}

\section*{Acknowledgements}
The research of Marc Hallin and Davy Paindaveine is supported by the
IAP research network grant \mbox{nr.} P7/06 of the Belgian government
(Belgian Science Policy). Davy Paindaveine moreover benefits from an
A.R.C. contract of the Communaut\'{e} Fran\c{c}aise de Belgique. The
research of Marc Hallin and Zudi Lu is supported by a Discovery Grant,
and Zudi Lu benefits from a Future Fellowship Grant, both of the
Australian Research Council. Miroslav \v Siman acknowledges the support
of Project 1M06047 of the Ministry of Education, Youth and Sports of
the Czech Republic. Marc Hallin and Davy Paindaveine are members of
ECORE, the association between CORE and ECARES. Marc Hallin is also
member of the Acad\' emie Royale de Belgique, and a fellow of CentER (Tilburg University).

%\begin{supplement}%[id=suppA]
%\slink[doi]{10.3150/00-BEJXXXXSUPP} %[doi,text={...}] - jei reikia
%suskaldyti doi
%\sfilename{BEJ000\_supp.pdf}
%\end{supplement}

% imsref loaded by jurgita.kaciuliene, 2014-04-16 15:44:45
%

\printhistory

\end{document}